\documentclass[10pt]{amsart}
\usepackage{amsfonts}
\usepackage{color}

\definecolor{c20}{rgb}{0.,0.7,0.}
\definecolor{c30}{rgb}{0.,0.,1.}
\definecolor{c40}{rgb}{1,0.1,0.7}
\definecolor{c50}{rgb}{1,0,0}
\definecolor{c60}{rgb}{1,0.9,0.1}

\def\EHH#1{\textcolor{c30}{#1}}
\def\EHH#1{#1}

\def\aE#1{\textcolor{c20}{#1}}

\def\EE#1{\textcolor{c40}{#1}}

\def\cL#1{\textcolor{c50}{#1}}
\def\cL#1{#1}

\newcommand{\tm}[1]{{\textcolor{green}{#1}}}
\def\tm#1{#1}

\def\cxL#1{\textcolor{c50}{#1}}
\def\cxL#1{#1}
\def\eW#1{\textcolor{c30}{#1}}
\def\eW#1{#1}
\def\EE#1{#1}
\def\aE#1{#1}
\newcommand{\kb}[1]{\boldsymbol{#1}}
\newcommand{\vk}[1]{\kb{#1}}

\newcommand{\ve}{\varepsilon}
\def\fracl#1#2{\biggr(\frac{#1}{#2} \biggl) }

\newcommand{\abs}[1]{\left\lvert #1 \right\rvert}
\newcommand{\Abs}[1]{ \biggl \lvert #1 \biggr \rvert}

\newcommand{\E}[1]{\mathbb{E}\left\{#1\right\}}

\newcommand{\pk}[1]{\mathbb{P} \left \{#1 \right \} }

\newcommand{\R}{\mathbb{R}}

\newcommand{\N}{\mathbb{N}}
\newcommand{\inr}{\in \R}
\newcommand{\inn}{\in \N}

\newcommand{\limit}[1]{\lim_{#1 \to   \infty}}

\newcommand{\BQN}{\begin{eqnarray}}
\newcommand{\EQN}{\end{eqnarray}}
\newcommand{\BQNY}{\begin{eqnarray*}}
\newcommand{\EQNY}{\end{eqnarray*}}

\newcommand{\BS}{\begin{sat}}
\newcommand{\ES}{\end{sat}}
\newcommand{\BT}{\begin{theo}}
\newcommand{\ET}{\end{theo}}
\newcommand{\BL}{\begin{lem}}
\newcommand{\EL}{\end{lem}}
\newcommand{\BK}{\begin{korr}}
\newcommand{\EK}{\end{korr}}

\newcommand{\BD}{\begin{de}}
\newcommand{\ED}{\end{de}}

\newcommand{\BIT}{\begin{itemize}}
\newcommand{\EIT}{\end{itemize}}
\newcommand{\BDI}{\begin{description}}
\newcommand{\EDI}{\end{description}}

\newcommand{\BRM}{\begin{remarks}}
\newcommand{\ERM}{\end{remarks}}

\newcommand{\BEL}{\begin{lem}}
\newcommand{\EEL}{\end{lem}}

\newtheorem{theo}{Theorem}[section]
\newtheorem{sat}[theo]{Proposition}
\newtheorem{de}[theo]{Definition}
\newtheorem{lem}[theo]{Lemma}

\newtheorem{example}[theo]{Example}
\newtheorem{korr}[theo]{Corollary}
\newtheorem{remark}[theo]{Remark}
\newtheorem{remarks}[theo]{Remarks}
\newtheorem{prop}[theo]{Proposition}

\newcommand{\nelem}[1]{{Lemma \ref{#1}}}
\newcommand{\neprop}[1]{{Proposition \ref{#1}}}
\newcommand{\netheo}[1]{{Theorem \ref{#1}}}
\newcommand{\nekorr}[1]{{Corollary \ref{#1}}}

\newcommand{\prooftheo}[1]{ \textsc{\bf Proof of Theorem} \ref{#1}:}
\newcommand{\proofprop}[1]{\textsc{\bf Proof of Proposition} \ref{#1}:}
\newcommand{\prooflem}[1]{\textsc{\bf Proof of Lemma} \ref{#1}:}

\newcommand{\COM}[1]{}

\newcommand{\QED}{\hfill $\Box$}

\def\IF{\infty}

\newcommand{\expon}[1]{\exp\left(#1\right)}

\date{}

\def\LT{\left}
\def\RT{\right}

\def\equivdis{\stackrel{d}{=}}

\newcommand{\toas}{\stackrel{a.s.}{\to}}

\def\u{\vk{u}}
\def\DU{ \Delta(\vk{u})}
\def\DUR{ \Delta_{(r)}(\vk{u})}
\def\RUR{ \Theta_{(r)}( \vk{u})}
\def\tilde{\widetilde}
\def\AIJIK{\EHH{Q_{ij,ik}}}
\def\AIJLK{Q_{ij,lk}}
\def\AIJIKP{\EHH{Q^+_{ij,ik}}}
\def\AIJLKP{Q^+_{ij,lk}}

\def\AARIJ{\abs{\arcsin(\sigma^{(1)}_{ij,lk})- \arcsin( \sigma^{(0)}_{ij,lk})}}
\def\AARIJ{\AIJLK}
\def\Yr{\EHH{\widetilde Y_{r:n}(t)}}
\begin{document}

\title[Comparison Inequalities for Order Statistics]{Comparison Inequalities for Order Statistics of Gaussian Arrays}

\author{Krzysztof D\c{e}bicki}
\address{Krzysztof D\c{e}bicki, Mathematical Institute, University of Wroc\l aw, pl. Grunwaldzki 2/4, 50-384 Wroc\l aw, Poland}
\email{Krzysztof.Debicki@math.uni.wroc.pl}

\author{Enkelejd  Hashorva}
\address{Enkelejd Hashorva, Department of Actuarial Science, 
University of Lausanne,\\
UNIL-Dorigny, 1015 Lausanne, Switzerland
}
\email{Enkelejd.Hashorva@unil.ch}

\author{Lanpeng Ji}
\address{Lanpeng Ji, Department of Actuarial Science, 
University of Lausanne\\
UNIL-Dorigny, 1015 Lausanne, Switzerland
}
\email{Lanpeng.Ji@unil.ch}

\author{Chengxiu Ling}
\address{Chengxiu Ling, Department of Actuarial Science, 
University of Lausanne\\
UNIL-Dorigny, 1015 Lausanne, Switzerland
}
\email{Chengxiu.Ling@unil.ch}

\bigskip

\date{\today}
 \maketitle

{\bf Abstract:}  Normal comparison lemma and Slepian's inequality are essential tools in
the study of Gaussian processes. In this paper we extend normal comparison lemma and
derive various related comparison inequalities including Slepian's inequality for order
statistics of two Gaussian arrays.The derived results can be applied 
in numerous problems related to the study of the supremum of
order statistics of Gaussian processes.
In order to illustrate the range of possible applications,
we analyze the lower tail behaviour of order statistics
of self-similar Gaussian processes and  derive mixed Gumbel limit theorems.
\\
{\bf Key Words:} comparison inequality; order statistics process;
Slepian's inequality; mixed Gumbel limit theorem; lower tail probability; self-similar Gaussian process; fractional Brownian motion.

\section{Introduction}\label{sec1}
The normal comparison inequality is crucial for the study of extremes of Gaussian processes, chi-processes and Gaussian random fields; see, e.g., \cite{Berman82, Berman92, leadbetter1983extremes, LiShao2004, Pit96}. It has been shown \tm{to be} valuable in many other fields of mathematics, \tm{such as, for instance, certain 
problems in number theory;} see,  e.g., \cite{Harper2, Harper1}.
  In the simpler framework of two $d$-dimensional  Gaussian distributions $\Phi_{\cxL{\Sigma}^{(1)}}$ and $\Phi_{\cxL{\Sigma}^{(0)}}$ with $N(0,1)$ marginal distributions,  the normal comparison inequality gives
bounds for the difference
$$\DU:= \Phi_{\EHH{\Sigma}^{(1)}}(\u)-\Phi_{\EHH{\Sigma}^{(0)}}(\u),\quad \forall \u=(u_1, \ldots, u_d)\inr^d$$
{by} a function of the covariance matrices
${\EHH{\Sigma}^{(k)}}=(\sigma^{(k)}_{ij})_{d\times d},k=0,1$. \EHH{As mentioned in} \cite{Kratz06}, \cxL{the derivation of} the bounds for $\Delta(\u)$, 
\cxL{by}  Slepian \cite{Slepian62}, Berman \cite{Berman64,Berman92}, Cram\'{e}r
\tm{\cite{Cramer67}}, \EHH{Bickel and Rosenblatt \cite{BickelRA}} and Piterbarg \cite{Pit72, Pit96}\cL,
\tm{relies} strongly on Plackett's partial differential equation; see \cite{Plackett}.  
 The most elaborated version of  the normal comparison  inequality is due to  Li and Shao  \cite{LiShao02}. Specifically, Theorem 2.1  therein shows that
\BQNY
 \DU \leq \frac{1}{2 \pi }\sum_{1\le i< j\le d}
\left(\arcsin(\EHH{\sigma}^{(1)}_{ij})- \arcsin( \EHH{\sigma}^{(0)}_{ij})\right)_+
\expon{-\frac{u^2_{i}+u^2_{j}}{2(1+\rho_{ij})}}, \quad \forall \u\inr^d,
\EQNY
where
$ \rho_{ij}:=\max(|\EHH{\sigma}^{(0)}_{ij}|,|\EHH{\sigma}^{(1)}_{ij}|)$ and $x_+=\max(x,0)$.
Clearly, if $\EHH{\sigma}^{(0)}_{ij} \ge \EHH{\sigma}^{(1)}_{ij}, 1\le i,j\le d$, then $\Phi_{\EHH{\Sigma}^{(1)}}(\u) \le \Phi_{\EHH{\Sigma}^{(0)}}(\u)$, which
is the well-known Slepian's inequality derived in \cite{Slepian62}.
\COM{
 In the meanwhile, under some technical conditions  Theorem 2.2 therein gives sharper bounds for
$ {\Phi_{\EHH{\Sigma}^{(1)}}(\u)}/{\Phi_{\EHH{\Sigma}^{(0)}}(\u)}$
when $u_i\ge0,1\le i\le d$, which is the main result of their paper. Later on,
following }
\EHH{Based on} the results of  \tm{Li and Shao \cite{LiShao02}, Yan \cite{Yan2009}} \EHH{showed that} for \eW{$\mathcal{N}$ an $N(0,1)$ random variable}
\BQN \label{ratio}
\quad 1\le  \frac{\Phi_{\EHH{\Sigma}^{(1)}}(\u)}{\Phi_{\EHH{\Sigma}^{(0)}}(\u)}\le  \exp\LT( \frac 1 {\sqrt{2 \pi}} \sum_{1\le i< j\le d} \frac{e^{-(u_i+u_j)^2/8}}{ \E{(\mathcal{N}+ (u_i+u_j)/2)_+ }}   \ln\LT(\frac{\pi-2\arcsin(\EHH{\sigma}^{(0)}_{ij})}{\pi-2\arcsin(\EHH{\sigma}^{(1)}_{ij})}\RT)\RT), \quad \vk{u} \in (0,\IF)^d
\EQN
provided that  $0\le \EHH{\sigma}^{(0)}_{ij}\le \EHH{\sigma}^{(1)}_{ij}\le 1$.
 Recent extensions of the normal comparison inequalities are presented in  \cite{DebickiHJL14, Harper2, Harper3, HashWeng, LuW2014, Monhor2013}.

In this paper, we are interested in the derivation of  comparison inequalities for
order statistics of Gaussian arrays, which are useful in several applications. To fix the notation, we denote by
$\mathcal{X}=(X_{ij})_{d\times n}$ and $\mathcal{Y}=(Y_{ij})_{d\times n}$
two random \tm{$d\times n$} arrays with $N(0,1)$ components
{(hereafter referred to as standard Gaussian arrays)}, and let
$\Sigma^{(1)}=(\sigma^{(1)}_{ij,lk})_{{dn \times dn}}$ and $\Sigma^{(0)}=(\sigma^{(0)}_{ij,lk})_{dn \times dn}$ be the covariance matrices of $\mathcal{X}$ and $\mathcal{Y}$, respectively, with
$ \sigma^{(1)}_{ij,lk} := \E{X_{ij}X_{lk}}$ and $\sigma^{(0)}_{ij,lk} := \E{Y_{ij}Y_{lk}}$. 
Furthermore, define $\vk{X}_{(r)}=(X_{1(r)}, \ldots, X_{d(r)}), 1\le r\le n$ to be the  $r$th order statistics vector generated by $\mathcal{X}$ as follows
\BQNY
X_{i(1)}= \min_{1\le j\le n} X_{ij} \le  \cdots \le X_{i(r)}\le\cdots \le \max_{1\le j\le n}X_{ij}= X_{i(n)},\quad 1\le i\le d.
\EQNY
Similarly, we write    $\vk{Y}_{(r)}
=(Y_{1(r)}, \ldots, Y_{d(r)})$ which is generated by $\mathcal{Y}$.

Our principal results, \netheo{T1} and \netheo{T3} \EHH{derive bounds} for the difference
\BQN\label{eq:Diff}
\DUR:= \pk{\vk{X}_{(r)} \le \u}- \pk{\vk{Y}_{(r)} \le \u}, \quad\u\in\R^d.
\EQN 
\EHH{A direct application of those bounds concerns the study of supreum of the $r$th order statistics process}
$\{X_{r:n}(t), t\ge0\}$ of $\{X_j(t), t\ge0\}, 1\le j\le n$ which are independent copies of a centered Gaussian process $\{X(t), t\ge0\}$. \EHH{More precisely, $X_{r:n}$ is defined by}
\BQNY
X_{n:n}(t){=\min_{1\le i \le n}X_i(t)}
 \le \cdots \le X_{r:n}(t)\le \cdots \le {\max_{1\le i \le n} X_i(t) =} X_{1:n}(t), \quad t\ge0.
\EQNY
\EHH{Below we call $X_{r:n}$ the $r$th order statistics process generated by $X$;   we} refer to \cite{DebickiHJL14, DebickiHJminima, DHJT15}  for the study of the extremes of order statistics processes.
\\
\EHH{With motivation from Theorem 3.1 in \cite{LiShao2004}, we apply the findings of \netheo{T1} to show that}
\BQN\label{fBm.r}
    \pk{ \sup_{t \in [0,1]} B_{r:n,\alpha}(t) \le x}= x^{ 2 p_{r:n,\alpha}/\alpha + o(1)},\quad x\downarrow 0
\EQN
holds with some non-negative constant $p_{r:n,\alpha}$, \EHH{where $B_{r:n,\alpha}$ is the $r$th order statistics process generated by
a fractional Brownian motion (fBm) $B_\alpha$ with Hurst index $\alpha/2 \in (0,1)$. }
\EHH{Moreover, if $B_{\alpha}^{(0)}$ is \tm{an fBm which is independent of}
$ B_{r:n,\alpha}$}, \EHH{then}
\BQN\label{fBm.rp}
    \pk{ \sup_{t \in [0,1]} \big(B_{r:n,\alpha}(t)-B_{\alpha}^{(0)}(t)\big) \le x}= x^{ 2 q_{r:n,\alpha}/\alpha + o(1)},\quad x\downarrow 0
\EQN
holds with some non-negative constant $q_{r:n,\alpha}$.
This result \tm{is} related to the problem of \tm{the capture time of a}
fractional Brownian pursuit;  see Theorem 4.1 in \cite{LiShao2004}.
\COM{Moreover, as an application of \netheo{T3} we derive, in \netheo{ThmB}, {several}
 limit theorems for the supremum of the order statistics processes generated by a stationary Gaussian process.}

\EHH{In \neprop{T2} we derive} bounds for \EHH{the ratio}
\BQNY
 \RUR:=\frac{\pk{\vk{X}_{(r)} \le \u } }{\pk{\vk{Y}_{(r)} \le \u }}, \quad \forall \u\in[0,\IF)^d,
\EQNY
which \tm{extend \eqref{ratio}}.
Relying on the findings of \EHH{Li and Shao}  \cite{LiShao02},
results for $\RUR$ can be \tm{applied for estimation of} $p_{r:n,\alpha}$ and $q_{r:n,\alpha}$ appearing in \eqref{fBm.r} and \eqref{fBm.rp}, respectively; \EHH{this topic will be investigated in a forthcoming paper. }

We organize this paper as \tm{follows.}
In Section \ref{Mainre} we \EHH{display}  our main results.
Section 3 is devoted to the study of the lower tail probability
of order statistics of self-similar Gaussian processes,
where  
\tm{Slepian-type} inequalities for order statistics processes are also derived.
We present the limit theorems for stationary order statistics processes in Section 4. Finally,
 all the proofs are relegated to Section \ref{sec4} \EHH{and Appendix.}

\section{ Main Results} \label{Mainre}

\tm{We begin this section with deriving some} 
sharp bounds for  $\DUR$ defined in \eqref{eq:Diff}\tm{, which go in line with}
Li and Shao's \cite{LiShao02} normal comparison inequality.
For notational simplicity we set below
$$\AIJLK:=\abs{\arcsin(\sigma^{(1)}_{ij,lk})- \arcsin( \sigma^{(0)}_{ij,lk})},  \quad
\AIJLKP:=( \arcsin(\sigma^{(1)}_{ij,lk})- \arcsin( \sigma^{(0)}_{ij,lk}))_+.$$

\BT\label{T1}
If
$\mathcal{X} $ and $\mathcal{Y}$ \EHH{are} two \tm{\cxL{standard} 
$d\times n$} 
Gaussian arrays, then for any $1\le r\le n$ \EHH{we have}

\BQN \label{ineq: abs}
\qquad\ \abs{\DUR}  \le \frac1{2\pi}\left( \underset{1\le j<k\le n}{\sum_{1\le i\le d}} \AIJIK \expon{-\frac{u_i^2}{1+\rho_{ij,ik}}}  + \underset{1\le   j,k\le n}{\sum_{1\le i<l \le d}} \AARIJ\expon{-\frac{u_i^2+u_l^2}{2(1+\rho_{ij,lk})}}\right),\  \forall \u\inr^d,
\EQN
where $\rho_{ij,lk}:=\max(|\sigma^{(0)}_{ij,lk}|,|\sigma^{(1)}_{ij,lk}|)$. \EHH{If further}
\BQN\label{cond: s-ind}
\sigma_{ij,ik}^{(1)} = \sigma_{ij,ik}^{(0)},\quad  1\le i\le d,\ 1\le j, k\le n,
\EQN
then
\BQN \label{ineq: without abs}
\DUR  \le  \frac 1{2\pi}\underset{1\le   j,k\le n}{\sum_{1\le i<l \le d}}
 \AIJLKP
\expon{-\frac{u_i^2+u_l^2}{2(1+\rho_{ij, lk})}},\quad \forall \u\inr^d.
\EQN
\ET

{ \remark \label{rem00}

a)  For $r=1$ and $r=n$ the claims in \eqref{ineq: abs} have been obtained in \cite{DebickiHJL14}, Lemma 4.1.
In fact, {using in addition similar arguments as in Theorem 1.2 in \cite{Pit96}}\cL, one can establish  for any  $[\boldsymbol a, \boldsymbol b]\subset [-\IF,\IF]^d$
the following comparison inequality
\BQNY
\lefteqn{\abs{\pk{\vk X_{(r)}\in[\boldsymbol a, \boldsymbol b]} - \pk{\vk Y_{(r)}\in[\boldsymbol a, \boldsymbol b]} }}\\
&&  \le  \frac{1}{\pi} \left( \underset{1\le j<k\le n}{\sum_{1\le i\le d}} \AIJIK
\expon{-\frac{u_i^2}{1+\rho_{ij,ik}}} + \underset{1\le   j,k\le n}{\sum_{1\le i<l \le d}} \AARIJ\expon{-\frac{u_i^2+u_l^2}{2(1+\rho_{ij,lk})}}\right) 
\EQNY
with  $ u_i=\min(|a_i|, |b_i|), 1\le i\le d.$\\
b) 
If  {$u_i\ge0, 1\le i\le d$} are all large enough, then 
\BQN\label{Rem:1} 
\DUR  \le  \frac{\cL1}{2\pi}\left(\underset{1\le j<k\le n}{\sum_{1\le i\le d}}
\AIJIKP 
\expon{-\frac{u_i^2}{1+\rho_{ij,ik}}}  +\underset{1\le   j,k\le n}{\sum_{1\le i<l \le d}}
\AIJLKP
\expon{-\frac{u_i^2+u_l^2}{2(1+\rho_{ij,lk})}}\right).
\EQN
{The proof of \eqref{Rem:1} is presented in the Appendix.} }

\EHH{A direct consequence} of   \netheo{T1} is the following  Slepian's inequality for the order statistics of Gaussian arrays.

\BK
If \eqref{cond: s-ind} is satisfied and further
$$\sigma_{ij,lk}^{(0)} \ge  \sigma_{ij,lk}^{(1)}, \ \ 1\le i< l \le d, 1\le j, k\le n$$
holds, then
$\DUR \cL\le 0$, or equivalently
\BQN\label{ineq: discrete}
\cxL{\mathbb P\Big\{}\cup_{i=1}^d\{X_{i(r)} > u_i\}\Big\} \ge \cxL{\mathbb P\Big\{}\cup_{i=1}^d\{Y_{i(r)} > u_i\}\Big\},\quad  \forall \u\inr^d.
\EQN

\EK

Note that the bounds in  \netheo{T1} do not include $r$, which {indicates} that
\tm{in some case they may be not sharp.}
\EHH{Below} we present \tm{a sharper result which holds} under the assumption
that the columns of both $\mathcal{X}$ and $\mathcal{Y}$ are mutually independent, i.e.,
\BQN\label{Assump.ind}
\sigma_{ij,lk}^{(\kappa)}= \sigma_{il}^{(\kappa)}\mathbb I\{j=k\}, \quad  1\le i,  l\le d, 1\le j, k\le n, \kappa=0,1,
\EQN
with some $\sigma_{il}^{(\kappa)},  1\le i,  l\le d, \kappa=0,1$, where $\mathbb I\{\cdot\}$ stands for the indicator function.
\EHH{This result is useful for establishing mixed Gumbel limit theorems,} see Section 4.

\EHH{In order to simplify the presentation, we shall define}
 $$
 c_{n,r}=\frac{n!}{ r!(n-r)! },\ \cL0\le r\le n,\ \ \ \rho_{il} =\max(\lvert \sigma_{il}^{(0)}\lvert, \lvert\sigma_{il}^{(1)}\lvert),\quad 1\le i,l\le d
 $$
 and
\BQNY 
 A_{il}^{(r)} = \int_{\sigma_{il}^{(0)}}^{\sigma_{il}^{(1)}}\frac{(1+\abs h)^{2(n-r)}}{(1-h^2)^{(n-r+1)/2}}\, dh,\quad 1\le i,l\le d,\ 1\le r\le n.
\EQNY

\BT\label{T3}
\EHH{Under the assumptions of \netheo{T1}, if further \eqref{Assump.ind} is satisfied}, then  for any $\vk u\in {(0,\IF)^d}$
\BQN\label{Ineq: sharpest}
\DUR \le \frac{n(c_{n-1, r-1})^2}{(2\pi)^{n-r+1}} u^{-2(n-r)}\sum_{1\le i<l\le d}\cxL{\big(}A_{il}^{(r)}\big)_+\expon{-\frac{(n-r+1)u^2}{1+\rho_{il}}}  
\EQN
and 
\BQNY
\abs{\DUR}
 \le   \frac{n(c_{n-1, r-1})^2}{(2\pi)^{n-r+1}} u^{-2(n-r)}\sum_{1\le i<l\le d}\abs{A_{il}^{(r)}}\expon{-\frac{(n-r+1)u^2}{1+\rho_{il}}}
\EQNY
hold with $u=\min_{1\le i\le d}u_i$. 
\ET

\tm{In the following proposition we derive  an upper bound for $\RUR$,  see related results in \cite{LiShao02, LuW2014, Yan2009}.}
\BS \label{T2}
\EHH{Under the assumptions of \netheo{T1}, if further  \eqref{cond: s-ind} holds \cL{and} } 
$\sigma_{ij,lk}^{(1)} \ge  \sigma_{ij,lk}^{(0)}\ge0$ for all $1\le i< l \le d, 1\le j, k\le n$, then for any $\vk u \in [0,\IF)^d$
\BQN \label{ineq: T2}
1\le \RUR \le  \exp\left( \eW{\frac 1 {\sqrt{2 \pi}}}\underset{1\le   j,k\le n}{\sum_{1\le i<l \le d}}
\frac{ C_{ij,lk} e^{-(u_i+u_l)^2/8}}{ \E{(\mathcal{N} + (u_i+u_l)/2)_+}} \right)
\EQN
 with $\mathcal{N}$ an $N(0,1)$ random variable and 
 $$C_{ij,lk}= \ln\fracl{\pi-2\arcsin(\sigma_{ij,lk}^{(0)})} {\pi-2\arcsin(\sigma_{ij,lk}^{(1)})},\ \ \   1\le i\cL< l\le d, 1\le j, k\le n.$$ 
\ES

\COM{
 {\remark \label{Rem:0}  It follows from the proof of \neprop{T2}, combined with corresponding arguments in \cite{LiShao02, Yan2009}, that the upper bound of
$\RUR$ also holds with  $H\big((u_i+u_l)/2\big)$ and/or $C_{ij,lk}$ replaced by
$\expon{\frac{u_i^2 +u_l^2}{2(1+ \rho_{ij,lk})}}$ and/or $ 2(\arcsin(\sigma_{ij,lk}^{(1)}) - \arcsin(\sigma_{ij,lk}^{(0)}))/\pi$, respectively.
 }
 }



\section{Lower Tail Probabilities of Order Statistics Processes}\label{sec3}

\EHH{The seminal contributions \cite{LiShao2004, LiShao2005} show that the investigation of}
the lower tail probability of Gaussian processes is of special interest in many applied fields,
for example, in the study of real zeros of random polynomials, in the study of fractional
Brownian pursuit, and in the study of \cxL{the} first\cxL{-}passage time for \cxL{the} Slepian process.
In this section, we aim at \EHH{extending} some results in \cite{LiShao2004, LiShao2005},
\tm{by considering} order statistics processes instead of Gaussian processes.

Our first result is concerned with Slepian's inequality for order statistics processes.
\EHH{In the following $X,Y,Z$ are three independent mean-zero Gausian processes with \cxL{almost surely (a.s.)} continuous sample paths.
 In accordance with our notation above
$X_{r:n}, Y_{r:n}$ and $Z_{r:n}$ are the corresponding $r$th order statistics processes.\\
Below,  we shall denote by
$\sigma_X(s,t) = \E{X(s)X(t)}$ and $\sigma_Y(s,t) = \E{Y(s)Y(t)}$ the  covariance function\cxL{s} of $X$ and $Y$\cxL{, respectively}.}

\begin{prop}\label{Prop1}
\COM{Let $\{X(t),t\ge0\}, \{Y(t), t\ge0\}$ be two  centered Gaussian processes with almost surely (a.s.) continuous sample paths, and let $\{Z(t), t\ge0\}$ be another one which is independent of them.} 
If for all $s,t\ge0$
\BQNY
\sigma_X(t, t)=\sigma_Y(t, t) \quad {\rm and}\ \ \sigma_X(s, t)\le \sigma_Y(s, t),
\EQNY
then for any $c\ge0,$ $T>0$ and  $ u\in\R$ \EHH{we have}
\BQN\label{Slep:Order1}
\pk{\sup_{t\in[0, T]}\big( X_{r:n}(t) +cZ(t)\big)>u} \; \ge \; \pk{\sup_{t\in[0, T]}\big( Y_{r:n}(t)+cZ(t)\big)>u}
\EQN
{and
\BQN\label{Slep:Order2}
\pk{\sup_{t\in[0, T]}\big( Z_{r:n}(t) +cX(t)\big)>u} \; \ge \; \pk{\sup_{t\in[0, T]}\big( Z_{r:n}(t)+cY(t)\big) >u}.
\EQN
}
\end{prop}
We shall display the proof of this proposition in Appendix. 


{\remark 
\COM{b) A natural generalization of Slepian's inequality  is the {claim in \eqref{Slep:Order1} with $c=0$.} Here in order to deal with the fractional Brownian pursuit  problem we allow $c>0$, which also means that Slepian's inequality holds for the order statistics processes of {(linearly) dependent} Gaussian processes. 
}
 \EHH{Similarly to  \cite{LiShao2004}, as \cL{an} immediate consequence of  Proposition \ref{Prop1}}
%
%
for any $c\ge0$  and $x\in\R$ we obtain that
\BQN\label{pr}
\quad p_r(x):=\lim_{T\to\IF}\frac1T\ln \pk{\sup_{0\le t\le T}\big( X_{r:n}(t) +cZ(t)\big)  \le x} =\sup_{T>0}\frac1T\ln \pk{\sup_{0\le t\le T}\big( X_{r:n}(t) +cZ(t)\big)\le x}
\EQN
exists and $p_r(x), x\inr$ is left-continuous, provided that $\sigma_X(0,t) \ge 0$ and $\sigma_Z(0,t) \ge 0$ for all $t\ge0$.
}

\EHH{ Next, we present} \tm{the main result of this section, which
gives a lower tail probability for} order statistics processes.

\tm{Let } $\{X(t), t\ge0\}$ be a centered  self-similar Gaussian process with   a.s. continuous sample paths and index $\alpha/2$ for some $\alpha>0$, i.e.,
\BQNY \label{def:selfsimilar}
X(0)=0, \quad \E{X^2(1)}=1, \quad \{X(\lambda t), t\ge0\}\stackrel d{=}\{\lambda^{\alpha/2} X(t), t\ge0\}, \quad \forall \lambda>0,
\EQNY
where $\stackrel d{=}$ denotes {equality of the (finite-dimensional) distribution functions.} 
It is well-known that by Lamperti's transformation a dual stationary Gaussian process $\{X^*(t), t\ge0\}$ can be defined as
$$
X^*(t)=e^{-\alpha t/2} X(e^t),\ \ t\ge0.
$$
\BT\label{ThmC}
Let $\{X(t),t\ge0\}$ and $ \{Z(t),t\ge0\}$ be two independent centered self-similar Gaussian processes with  continuous sample paths and common \EHH{self-similarity} index $\alpha/2\in \cL{(0, \IF)}$. 
Suppose that $\sigma_X(s,t)\ge0$ and  $\sigma_Z(s,t)\ge 0$ for all $s,t\ge0,$ and both $\rho(t):=  \E{X^*(t)X^*(0)}, 
\widetilde{\rho}(t):=  \E{Z^*(t)Z^*(0)}, t\ge0$ are  non-increasing. If further
for any $h\in(0,\IF)$ \cxL{and} $\theta\in(0,1)$
\BQN\label{continuity}
a_{h,\theta}^2= \inf_{0<t\le h}\frac{\rho(\theta t)-\rho(t)}{1-\rho(t)}>0,
\EQN
then for any $c\ge0$
\BQN\label{Extend_fBm}
\pk{\sup_{0\le t\le 1}\big(X_{r:n}(t)+cZ(t)\big)\le x }=x^{2\mathrm c_{r:n,\alpha}/\alpha+o(1)},\quad  x\downarrow 0,
\EQN
\tm{where
$$\mathrm c_{r:n,\alpha}:=-\lim_{T\to\IF}\frac1T\ln \pk{\sup_{0\le t\le T}  \big(X^*_{r:n}(t) + cZ^*(t)\big)  \le 0}$$
is the Li-Shao type constant.}
\ET

{\remark  As discussed in \cite{LiShao2005} two examples of $\{X(t),t\ge0\}$ that satisfy all conditions of \netheo{ThmC} are the standard fBm  $B_\alpha$ and the \EHH{centered}  Gaussian process $\{X_\beta(t), t\ge0\}, \beta>0$ with
$$ \E{ X_\beta(s) X_\beta(t)}= \frac{2^\beta (st)^{(1+\beta)/2}}{(s+t)^\beta}, \quad s,t> 0.$$
\EHH{Moreover}, we have that \eqref{fBm.r} holds with $p_{r:n,\alpha}$ given by
$$
p_{r:n,\alpha}=-\lim_{T\to\IF}\frac1T\ln \pk{\sup_{0\le t\le T} B_{r:n,\alpha}(e^t)   \le 0}. 
$$
}

\COM{\example Another class of Gaussian processes are $\{X_\beta(t), t\ge0\}, \beta>0$ with covariance function
$$ \E{ X_\beta(s) X_\beta(t)}= \frac{2^\beta (st)^{(1+\beta)/2}}{(s+t)^\beta}, \quad s,t> 0.$$
We see that the conditions of \netheo{ThmC} with $\alpha=1$ hold for this family of Gaussian processes. A larger class of
Gaussian processes {$X_{\alpha,\beta}$} is studied in \cite{AurzadaJTE}. From Proposition 3.3 in \cite{LiShao2004}, the lower tail probability of $X_\beta$ with $\beta=1$
is derived, and it is related to the probability that a random polynomial has no real zero. If we consider the order statistics processes $X_{r:n,\beta}$ generated by $X_\beta$, then in view of \netheo{ThmC}, we have
\BQNY
\pk{ \sup_{t\in [0,1]}X_{r:n,\beta}(t) \le x}= x^{2\mathrm c_{r,\beta}+ o(1)}, \quad x\downarrow 0,
\EQNY
where $\mathrm c_{r,\beta}$ is some non-negative constant. }

To this end, we discuss a modification of the fractional Brownian pursuit problem considered
in \cite{LiShao02, LiShao2004}. Let $B_{\alpha}^{(k)}, 0 \le k\le n$ be independent
standard fBm\cxL{s}, and \tm{define} 
\BQN\label{eq:tau}
\tau_{r:n,\alpha}=\inf\{t\ge0: B_{r:n,\alpha}(t)-1=B_{\alpha}^{(0)}(t)\}, 
\EQN
 where $\{ B_{r:n,\alpha}(t) ,t\ge0\}$ is the $r$th order statistics process of $B_{\alpha}^{(k)}, 1 \le k\le n$. Then $\tau_{r:n,\alpha}$ can be viewed as a capture time in a random pursuit setting. Assume that a fractional Brownian prisoner escapes, running along the path of $B_{\alpha}^{(0)}$. In his pursuit, there are $n$ independent fractional Brownian policemen running along the paths of $B_{\alpha}^{(k)}, 1 \le k\le n$, respectively. At the outset, the prisoner is ahead of the policemen by 1 unit of distance. Then, $\tau_{r:n,\alpha}$ represents the capture time when at least $r$ policemen catches the prisoner. As shown in the aforementioned papers
 the study of  the
capture time of the fractional Brownian pursuit is related to the analysis of the lower tail probability of order statistics process
since
\BQNY
\pk{\tau_{r:n,\alpha}>s}=\pk{\sup_{t\in [0,1]}\Bigl(B_{r:n,\alpha}(t)-B_{\alpha}^{(0)}(t) \Bigr) \le s^{-\alpha/2}  }.
\EQNY
As an application of \netheo{ThmC} we have that \eqref{fBm.rp} holds with
$$ q_{r:n,\alpha}=-\lim_{T\to\IF}\frac1T\ln \pk{\sup_{0\le t\le T}  {\big(B_{r:n,\alpha}(e^t)- B_{\alpha}^{(0)}(e^t)\big)}  \le 0},
$$
\EHH{which leads} to the following result.
\BK If $\tau_{r:n,\alpha}$ is defined as in \eqref{eq:tau}, \EHH{then we have}
\BQNY
 \pk{\tau_{r:n,\alpha}>\cxL{s}}  = \cxL{s}^{-  q_{r:n,\alpha} + o(1)}, \quad   \cxL{s}\to\IF.
\EQNY
\EK

\section{Limit Theorems of Stationary Order Statistics Processes}\label{sec3.2}
In this section \tm{we suppose that}
$\{X_{r:n}(t),t\ge0\}$ to be the $r$th order statistics process generated by a centered stationary Gaussian process $\{X(t),t\ge0\}$ with \cxL{a.s.} continuous sample paths, \EHH{\cL{unit variance} 
} and \EHH{correlation} function $\rho(\cdot)$ satisfying
\BQN\label{corrr}
\rho(t) = 1- \abs t^\alpha + o(\abs t^\alpha),
\ t\to0 \  {\rm \ for \ some\ }  \alpha\in(0, 2] \ \ \text{and} \  \ \rho(t) < 1,\ \forall t\not=0. 
\EQN

From Theorem 1.1 in \cite{DebickiHJL14} or  Theorem 2.2  in \cite{DHJT15} for any $T>0$
\EHH{we have}
\BQN\label{Asym: finite-time}
\pk{ \sup_{t\in [0,T]}
X_{r:n}(t) > u}= T \mathcal A_{r,\alpha} c_{n,r} (2\pi)^{-\frac r2} u^{\frac2\alpha -r} \expon{-\frac{ru^2}2}(1+o(1)), \quad
u\to \IF,
\EQN
where   $\mathcal A_{r,\alpha}\in(0,\IF)$  is a positive constant. 
As a continuation of \cite{DebickiHJL14} we  establish below a limit theorem  $X_{r:n}$. 
\BT\label{ThmB}
Let $\{X_{r:n}(t),t\ge0\}$ be the $r$th order statistics process generated by $X$ as above. Suppose that \eqref{corrr} holds and 
 $\lim_{t\to\IF} \rho(t)\ln t=\gamma \in[0, \IF]$.\\
a) If $\gamma=0$, then
\BQNY
 \limit{T}\sup_{x\in \R} \Abs{ \pk{{a_{r, T}} \Big(\sup_{ t\in [0, T]} X_{r:n}(t)- {b_{r, T}}\Big) \le x} -
\expon{ - e^{-x}}} =0,
\EQNY
where, with $D:= (r/2)^{r/2 -1/ \alpha}c_{n,r}\mathcal A_{r,\alpha}(2\pi)^{-r/2}$
\BQN \label{Normalized constant}
 a_{r, T} = \sqrt{2r\ln T}, \quad b_{r, T} =\sqrt{(2/r)\ln T} + \frac 1{\sqrt{2r\ln T}} \left(
\left(\frac 1\alpha -\frac r2\right)\ln\ln T + \ln D\right), \quad T>e.
 \EQN
 {b) If $\gamma=\IF$},  and $\alpha\in(0,1]$, $\rho(t)$ is convex for $t\ge0$ with  $\lim_{t\to\IF}\rho(t)=0$  and further  $\rho(t)\ln t$ is monotone for large $t$, then \cxL{with $\Phi(\cdot)$ the df of an $N(0,1)$ random variable,} 
\BQNY
 \limit{T}\sup_{x\in \R} \Abs{ \pk{ \frac1{\sqrt{\rho(T)}} \Big(\sup_{ t\in [0, T]} X_{r:n}(t)- \sqrt{1-\rho(T)}b_{r, T}\Big) \le x} -
\Phi(x)} =0.
\EQNY
c)  If $\gamma\in(0,\IF)$,  
\cxL{then, with $W$ an  $N(0,1)$ random variable}
\BQNY
 \limit{T}\sup_{x\in \R} \Abs{ \pk{ a_{r, T} \Big(\sup_{ t\in [0, T]} X_{r:n}(t)- b_{r, T}\Big) \le x} -
\E{\expon{ - e^{-(x+\gamma -\sqrt{2\gamma r}W\cL)}}}} =0.
\EQNY
\ET
\tm{The proof of Theorem \ref{ThmB} is presented in the Appendix.}

\COM{{\remark \label{rem:ThmB}
Using the crucial normal comparison theorem for Gaussian order statistics (\netheo{T3}), following  similar arguments of Theorems 2.3, 3.1 in \cite{MittalY1975} to consider
alternatively the strong dependent conditions: $\lim_{t\to\IF}\rho(t) \ln t=\gamma\in(0, \IF]$  and keep other conditions in \netheo{ThmB}, we have\\
a) If $\gamma<\IF$, then with $W$ an N(0,1) random variable
\BQNY
 \limit{T}\sup_{x\in \R} \Abs{ \pk{ a_{r, T} \Big(\sup_{ t\in [0, T]} X_{r:n}(t)- b_{r, T}\Big) \le x} -
\E{\expon{ - e^{-(x+\gamma -\sqrt{2\gamma r}W\cL)}}}} =0.
\EQNY
b) If $\gamma=\IF$ and $\alpha\in(0,1]$ and further that $\rho(t)$ is convex for $t\ge0$ and  $\lim_{t\to\IF}\rho(t)=0$, and that $\rho(t)\ln t$ is monotone for large $t$, then, 
\cxL{with $\Phi(\cdot)$ the distribution function of an $N(0,1)$ random variable,}
\BQNY
 \limit{T}\sup_{x\in \R} \Abs{ \pk{ \frac1{\sqrt{\rho(T)}} \Big(\sup_{ t\in [0, T]} X_{r:n}(t)- \sqrt{1-\rho(T)}b_{r, T}\Big) \le x} -
\Phi(x)} =0.
\EQNY
Here $a_{r, T}$ and $b_{r, T}$ are given by \eqref{Normalized constant}. See also Theorems 2.1, 2.2 in \cite{TanHP2012} for more related discussions.
}
}

\COM{\remark We see that for the stationary order statistics processes (other than Gaussian processes), the Berman's condition is still acting as in the enough {and necessary} to ensure Gumbel limit theorems {as well as the mixing Gumbel limit theorem under strong dependent conditions}. The refined comparison inequality stated in \netheo{T3} is a powerful tool in establishing \netheo{ThmB} above.}

\section{Proofs}\label{sec4}
\COM{Section \ref{proof.T1} contains the proof of \netheo{T1} followed by the proofs of \neprop{T2} and \netheo{T3} whereas the order statistics $X_{i(r)}, Y_{i(r)}$ generated respectively by two normal random matrices $(X_{ij})_{d\times n}, (Y_{ij})_{d\times n}$. Finally, in Section \ref{proof.Prop1B} we give the proofs of \neprop{Prop1} and \netheo{ThmB} whereas we consider the order statistics process $\{X_{r:n}(t), t\ge0\}$ generated by $n$ centered Gaussian processes.}
\COM{This section is devoted to the proofs of all the results.
The proofs are based  on several ideas from \cite{ DebickiHJL14, LuW2014,Yan2009}. However,
it turns out that dealing with  order statistics (processes) is much more difficult. 
}
\EHH{Hereafter}, we write $\stackrel{d}{=}$ for \EHH{equality of the distribution functions.}
A vector $\vk z=(z_1,z_2, \ldots, z_{dn})$ will also be denoted by
$$
\vk{z}=(\vk{z}_1,\ldots,\vk z_d),\ \text{with\ } \vk z_i=(z_{i1},z_{i2},\ldots,z_{in}),\ \ 1\le i\le d,
$$
where $z_{ij}=z_{(i-1)n+j}, 1\le i\le d, 1\le j\le n$. Note that for any $p=(i-1)n+j, q=(l-1)n+k, 1\le i,l\le d, 1\le j,k\le n$
$$
\{p<q\}=\{i<l,\ \text{or\ } i=l \ \text{and\ } j<k  \}.
$$
Furthermore, for any $\vk{z}\in \R^n$ we denote
$$
\frac{d\vk z}{d z_i}=dz_1dz_2\cdots dz_{i-1}dz_{i+1}\cdots dz_n,\quad \cxL{1\le i \le n},
$$
and for $\cxL{1\le} i< j\cxL{\le n}$
$$
\frac{d\vk z}{d z_i d z_j}=dz_1dz_2\cdots dz_{i-1}dz_{i+1}\cdots  dz_{j-1}dz_{j+1}\cdots dz_n.
$$

\prooftheo{T1} {We shall first establish  \eqref{ineq: abs} by considering $r=1,$ $r=2$ and $2< r\le n$ separately.} \\
\underline{Case $r=1$}.
Note that $\mathcal{X} \equivdis -\mathcal{X}$ for the standard  Gaussian array $\mathcal{X}$. 
It follows from  Theorem 2.1 in \cite{LuW2014}  that
\begin{align*}
\lefteqn{\abs{\Delta_{(1)}(\vk u)}=
\Abs{\cxL{\mathbb P\Big\{}\cup_{i=1}^d\cap_{j=1}^n\{-Y_{ij}<-u_i\}\Big\} - \cxL{\mathbb P\Big\{}\cup_{i=1}^d\cap_{j=1}^n\{-X_{ij}<-u_i\}\Big\}}}
\notag\\
&\le  \frac1{2\pi}\sum_{(i-1)n+j<(l-1)n+k}  \AARIJ \expon{-\frac{u_i^2+u_l^2}{2(1+\rho_{ij,lk})}}
\end{align*}
establishing the proof of $r=1$.\\

Next, by a standard approximation procedure  we may assume that both $\Sigma^{(1)}$ and $\Sigma^{(0)}$ are positive definite.
 Let further $ \mathcal  Z = (Z_{ij})_{d\times n}$ be a 
 standard Gaussian array  with
covariance matrix $ {\Gamma^h }= h\Sigma^{(1)} +(1-h)\Sigma^{(0)} =(\delta_{ij,lk}^h)_{dn \times dn}$, where
by our notation $\delta_{ij,lk}^h = \E{Z_{ij}Z_{lk}}$.
Clearly,  $\Gamma^h$ is also positive definite for any $h\in[0,1]$.
Denote below by $g_h(\vk z)$ the probability density function (pdf) of $\mathcal Z$. It is known that (see, e.g., \cite{leadbetter1983extremes}, p.\,82, or \cite{LuW2014}) 
\BQN\label{differentiate h}
\frac{\partial g_h(\vk z)}{ \partial \delta_{ij,lk}^h}=\frac{\partial^2 g_h(\vk z)}{ \partial z_{ij}\partial z_{lk}},\quad 1\le i, l\le d, 1\le j, k\le n.
\EQN

\underline{Case $r=2$}.  Hereafter, we write $\boldsymbol \lambda=-\vk u$ and set
\BQN\label{def: Q}
Q(\mathcal Z; \Gamma^h) = \pk{\vk Z_{(n-1)}>\boldsymbol\lambda} = \int_{\cap_{i=1}^d\cup_{j, j'=1}^n\{z_{ij} >\lambda_i, z_{ij'}>\lambda_i\}} g_h(\vk z)\, d\vk z.
\EQN
Since $\vk X_{(2)}\stackrel{d}{=}- \vk X_{(n-1)}$ we have
\BQN\label{Decomposition I}
\Delta_{(2)}(\vk u)= Q(\mathcal Z; \Gamma^1)- Q(\mathcal Z; \Gamma^0) = \int_0^1\frac{\partial Q(\mathcal Z; \Gamma^h)}{\partial h}\, dh.
\EQN
Note that the quantities $Q(\mathcal Z; \Gamma^h)$ and $g_h(\vk z)$ depend  on $h$ only through the entries $\delta_{ij,lk}^h$ of $\Gamma^h$. 
Hence we have by \eqref{differentiate h}
\BQN\label{Decomposition II}
{\frac{\partial Q}{\partial h}(\mathcal Z; \Gamma^h) }&=&  {\sum_{(i-1)n+j < (l-1)n +k}\frac{\partial Q(\mathcal Z; \Gamma^h)}{\partial \delta_{ij,lk}^h} \frac{\partial \delta_{ij,lk}^h}{\partial h}}  \notag \\
& = &\sum_{(i-1)n+j<(l-1)n+k}(\sigma_{ij,lk}^{(1)} -\sigma_{ij,lk}^{(0)}) {E_{il}(j,k)},
\EQN
with
$$
E_{il}(j,k):=\int_{\cap_{s=1}^d\cup_{t, t'=1}^n\{z_{st} >\lambda_s, z_{st'}>\lambda_s\}} \frac{\partial^2 g_h(\vk z)}{\partial z_{ij}\partial z_{lk}}\, d\vk z, \quad  (i-1)n+j<(l-1)n+k.
$$
Next, in order to establish \eqref{ineq: abs} we shall show that
\BQN\label{Decomposition IIIE}
\abs{E_{il}(j,k)} \le \varphi(\lambda_i, \lambda_l; \delta_{ij,lk}^h),\quad  (i-1)n+j<(l-1)n+k,
\EQN
where $\varphi(\cdot, \cdot; \delta_{ij,lk}^h)$ is the pdf of $(Z_{ij}, Z_{lk})$,
given by
\BQNY
\varphi(x,y; \delta_{ij,lk}^h) = \frac1{2\pi\sqrt{1-(\delta_{ij,lk}^h)^2}}
\expon{-\frac{x^2-2\delta_{ij,lk}^h x y +y^2}{2\tm{\left(1-(\delta_{ij,lk}^h)^2\right)}}},\quad x, y\in\R. 
\EQNY
We consider below   two sub-cases: a)  $i=l$ and b)  $i<l$.
\\
\underline{a) Proof of \eqref{Decomposition IIIE} for $i=l$}. Letting
$A_i\rq{}=\cap_{s=1,s\neq i}^d\cup_{t, t'=1}^n\{z_{st} >\lambda_s, z_{st'}>\lambda_s\}, A_i=\cup_{t,t\rq{}=1}^n\{z_{it}>\lambda_i, z_{it\rq{}}>\lambda_i\},
$ we rewrite $E_{ii}(j,k)$ as
\BQN\label{Eq: reexp}
E_{ii}(j,k) =\int_{A_i\rq{}}\int_{A_i}
\frac{\partial^2 g_h(\vk z)}{\partial z_{ij}\partial z_{ik}}\, d\vk z,\quad 1\le i\le d,\ 1\le j< k\le n.
\EQN
Next, we decompose the integral region $A_i$ according to
$a_1)\  z_{ij}>\lambda_i, z_{ik}>\lambda_i$; $a_2)\  z_{ij}>\lambda_i, z_{ik}\le\lambda_i$; $a_3)\  z_{ij}\le\lambda_i, z_{ik}>\lambda_i$; and
$a_4)\  z_{ij}\le\lambda_i, z_{ik}\le\lambda_i$. 

For case $a_1)$ we have
\BQN\label{Eq: a1}
\int_{A_i\cap\{z_{ij}>\lambda_i, z_{ik}>\lambda_i\}}
\frac{\partial^2 g_h(\vk z)}{\partial z_{ij}\partial z_{ik}}\, d\vk z_i
=\int_{\R^{n-2}}g_h(z_{ij}=\lambda_i, z_{ik}=\lambda_i)\, \frac{d\vk z_i}{dz_{ij} dz_{ik}},
\EQN
where $g_h(z_{ij}=\lambda_i, z_{ik}=\lambda_i)$ denotes a function of $dn-2$ variables formed from $g_h(\vk z)$ by putting $z_{ij}=\lambda_i, z_{ik}=\lambda_i$. 
Similarly, for cases $a_2)$ and $a_3)$
\BQN\label{Eq: a23}
\nonumber \lefteqn{\int_{A_i\cap\{z_{ij}>\lambda_i, z_{ik}\le\lambda_i\}}
\frac{\partial^2 g_h(\vk z)}{\partial z_{ij}\partial z_{ik}}\, d\vk z_i
=\int_{A_i\cap\{z_{ij}\le\lambda_i, z_{ik}>\lambda_i\}}
\frac{\partial^2 g_h(\vk z)}{\partial z_{ij}\partial z_{ik}}\, d\vk z_i } \\
&&=-
\int_{\cup_{t=1,t\neq j,k}^n\{z_{it}>\lambda_i\}}g_h(z_{ij}=\lambda_i, z_{ik}=\lambda_i)\, \frac{d\vk z_i}{dz_{ij} dz_{ik}}.
\EQN
Finally, for case $a_4)$
\BQNY
\int_{A_i\cap\{z_{ij}\le\lambda_i, z_{ik}\le\lambda_i\}}
\frac{\partial^2 g_h(\vk z)}{\partial z_{ij}\partial z_{ik}}\, d\vk z_i
=\int_{\cup_{t,t'=1,t,t'\neq j,k}^n\{z_{it}>\lambda_i,z_{it'}>\lambda_i\}}g_h(z_{ij}=\lambda_i, z_{ik}=\lambda_i)\, \frac{d\vk z_i}{dz_{ij} dz_{ik}}.
\EQNY
This together with 
\eqref{Eq: reexp}--\eqref{Eq: a23} yields
\BQN\label{Eq: Ee}
E_{ii}(j,k) &=&\int_{A_i\rq{}}\int_{\R^{n-2}-\cup_{t=1,t\neq j,k}^n\{z_{it}>\lambda_i\}}g_h(z_{ij}=\lambda_i, z_{ik}=\lambda_i)\, \frac{d\vk z}{dz_{ij} dz_{ik}}
\notag\\
&&\quad- \int_{A_i'}\int_{\cup_{t=1,t\neq j,k}^n\{z_{it}>\lambda_i\}-\cup_{t,t'=1,t,t'\neq j,k}^n\{z_{it}>\lambda_i,z_{it'}>\lambda_i\}}g_h(z_{ij}=\lambda_i, z_{ik}=\lambda_i)\, \frac{d\vk z}{dz_{ij} dz_{ik}} \notag\\
&=&
\lefteqn{\varphi(\lambda_i, \lambda_i; \delta_{ij,ik}^h)\left( \pk{(\cap_{s=1, s\neq i}^d\{Z_{s(n-1)} >\lambda_s\})\cap\{\vk Z_i''\in\{w_{i1}''=\IF\}\}\Big\lvert \{Z_{ij}=\lambda_i, Z_{ik} =\lambda_i\}}\right.} \notag \\
&&\quad\left. -\pk{(\cap_{s=1, s\neq i}^d\{Z_{s(n-1)} >\lambda_s\})\cap\{\vk Z_i''\in\{w_{i1}''\le n, w_{i2}''=\IF\}\}\Big\lvert \{Z_{ij}=\lambda_i, Z_{ik} =\lambda_i\}}\right),
\EQN
where 
$\vk Z_i''$ is the $(n-2)$-dimensional components of $\vk Z_i$ obtained by deleting $Z_{ij}$ and $Z_{ik}$, and $w_{i1}'', w_{i2}''$ are defined by (recall $\inf\{\emptyset\}=\IF$)
\BQN\label{w''}
w_{i1}''=\inf\{t: z_{it}>\lambda_i, t\neq j,k\}, \quad w_{i2}''=\inf\{t: z_{it}>\lambda_i,{t\neq j,k, t> w_{i1}''}\}.
\EQN
It follows from \eqref{Eq: Ee}  that 
\eqref{Decomposition IIIE} holds for $i=l$.
\\
\underline{b) Proof of \eqref{Decomposition IIIE} for $i< l$}.  With $A_{il}''=\cap_{s=1,s\neq i,l}^d\cup_{t, t'=1}^n\{z_{st} >\lambda_s, z_{st'}>\lambda_s\}$, we have (recall $A_i$ in \eqref{Eq: reexp})
\BQN\label{decomp. 1En}
E_{il}(j,k)=\int_{A_{il}''}\int_{A_l}\int_{A_i}
\frac{\partial^2 g_h(\vk z)}{\partial z_{ij}\partial z_{lk}}\, d\vk z.
\EQN
Next, we decompose the integral region $A_i$ according to $z_{ij}>\lambda_i$ and $z_{ij} \le\lambda_i$. We have
\BQNY
\lefteqn{ \int_{A_i\cap\{z_{ij}>\lambda_i\}}
\frac{\partial^2 g_h(\vk z)}{\partial z_{ij}\partial z_{lk}}\, d\vk z_i
+\int_{A_i\cap\{z_{ij}\le\lambda_i\}}
\frac{\partial^2 g_h(\vk z)}{\partial z_{ij}\partial z_{lk}}\, d\vk z_i
}\\
&&=-\int_{\cup_{t=1,t\neq j}^n\{z_{it}>\lambda_i\}-\cup_{t,t'=1,t,t'\neq j}^n\{z_{it}>\lambda_i, z_{it'}>\lambda_i\}}\frac{\partial g_h(z_{ij}=\lambda_i)}{\partial z_{lk}}\, \frac{d\vk z_i}{dz_{ij}}\\
&&=-\int_{\{w_{i1}'\le n, w_{i2}'=\IF\}}\frac{\partial g_h(z_{ij}=\lambda_i)}{\partial z_{lk}}\, \frac{d\vk z_i}{dz_{ij}},
\EQNY
where $w_{i1}', w_{i2}'$ are defined by (similar notation below for $w_{l1}', w_{l2}'$ with respect to $k$)
\BQN\label{w'}
w_{i1}'=\inf\{t: z_{it}>\lambda_i, t\neq j\}, \quad w_{i2}'=\inf\{t: z_{it}>\lambda_i, t\neq j, t> w_{i1}'\}.
\EQN
Using similar arguments for the integral with region $A_l$, we have by \eqref{decomp. 1En}
\BQN\label{Eq: En}
\notag E_{il}(j,k)&=&\int_{A_{il}''}\int_{\{ w_{i1}'\le n, w_{i2}'=\IF\}}\int_{\{w_{l1}'\le n, w_{l2}'=\IF\}} g_h(z_{ij}=\lambda_i, z_{lk}=\lambda_l)\, \frac{d\vk z}{dz_{ij} dz_{lk}}\\
&=&\varphi(\lambda_i, \lambda_l; \delta_{ij,lk}^h) \mathbb P\Bigl \{  \cap_{s=1, s\neq i,l}^d\{Z_{s(n-1)} >\lambda_s\}\cap (\vk Z_i'\in\{w_{i1}'\le n, w_{i2}'=\IF\}) \notag
\\ && \qquad\qquad \qquad\quad \cap(\vk Z_l'\in\{w_{l1}'\le n, w_{l2}'=\IF\})\Big\lvert \{Z_{ij}=\lambda_i, Z_{lk} =\lambda_l\}\Bigr\},
\EQN
where $\vk Z_i'$ and $\vk Z_l'$ are the $(n-1)$-dimensional components of $\vk Z_i$ and $\vk Z_l$ obtained by deleting $Z_{ij}$ and $Z_{lk}$, respectively. Consequently, by \eqref{Eq: Ee} and \eqref{Eq: En} the validity of \eqref{Decomposition IIIE} follows. Next, by combining
\eqref{Decomposition I}--\eqref{Decomposition IIIE}, the claim in \eqref{ineq: abs} for $r=2$ follows by the fact  that (see \cite{LiShao02})   
\BQN\label{Ineq: bivar norm}
\int_0^1 \varphi(\lambda_i,\lambda_l;\delta_{ij,lk}^h)\, dh \le
 \frac{\arcsin(\sigma_{ij,lk}^{(1)})-\arcsin(\sigma_{ij,lk}^{(0)})}{2\pi(\sigma_{ij,lk}^{(1)} - \sigma_{ij,lk}^{(0)})}\expon{-\frac{\lambda_i^2+\lambda_l^2}{2(1+\rho_{ij,lk})}}.
\EQN

\underline{Case $2<r\le n$}. 
Letting
$
\tilde Q(\mathcal Z;\Gamma^h)=\pk{\vk Z_{(n-r+1)}>\boldsymbol \lambda}$
we have
\begin{align}\label{Decomposition I(II)r}
\Delta_{(r)}( \vk u) 
 = \int_0^1 \, dh \left(\sum_{(i-1)n+j<(l-1)n+k} (\sigma_{ij,lk}^{(1)} -\sigma_{ij,lk}^{(0)}) \tilde E_{il}(j,k)\right),
\end{align}
where
$$
\tilde E_{il}(j,k):=\int_{\cap_{s=1}^d\cup_{t_1, \ldots, t_r=1}^n\{z_{st_1} >\lambda_s, \ldots, z_{st_r}>\lambda_s\}} \frac{\partial^2 g_h(\vk z)}{\partial z_{ij}\partial z_{lk}}\, d\vk z.
$$
With the aid of \eqref{Ineq: bivar norm}, it suffices to show that  
\BQN\label{Decomposition IIIET}
\abs{\tilde E_{il}(j,k)} \le \varphi(\lambda_i, \lambda_l; \delta_{ij,lk}^h),\quad  (i-1)n+j<(l-1)n+k.
\EQN
Similarly as above, two sub-cases : a)  $i=l$ and b)  $i<l$ need to be considered separately.

\underline{a) Proof of \eqref{Decomposition IIIET} for $i=l$}. Similarly to $E_{ii}(j,k)$, we rewrite $\tilde E_{ii}(j,k) $ as
\BQN\label{Eq: reexpp}
\tilde E_{ii}(j,k) =\int_{\tilde A_i\rq{}}\int_{\tilde A_i}
\frac{\partial^2 g_h(\vk z)}{\partial z_{ij}\partial z_{ik}}\, d\vk z,
\EQN
where
$\tilde A_i\rq{}=\cap_{s=1,s\neq i}^d\cup_{t_1, \ldots, t_r=1}^n\{z_{st_1} >\lambda_s, \ldots, z_{st_r}>\lambda_s\},\ \tilde A_i=\cup_{t_1, \ldots, t_r=1}^n\{z_{it_1} >\lambda_i, \ldots, z_{it_r}>\lambda_i\}.$\\
Next, we decompose the integral region $\tilde A_i$ according to the four cases
$a_1)$--$a_4)$ \tm{as} introduced for $A_i$ (see   the last two lines right above \eqref{Eq: a1}).

For case $a_1)$
\BQN\label{Eq: a1r}
\int_{\tilde A_i\cap\{z_{ij}>\lambda_i, z_{ik}>\lambda_i\}}
\frac{\partial^2 g_h(\vk z)}{\partial z_{ij}\partial z_{ik}}\, d\vk z_i
=\int_{\{w_{i,r-2}''\le n\}}g_h(z_{ij}=\lambda_i, z_{ik}=\lambda_i)\, \frac{d\vk z_i}{dz_{ij} dz_{ik}},
\EQN
where
$w_{i1}''$ is given by \eqref{w''} and (notation: $w_{i,t}'' = w_{it}''$)
\BQNY
w_{it}''=\inf\{t_0\le n: z_{it_0}>\lambda_i, t_0\neq j,k, t_0>w_{i,t-1}''\}, \quad 2\le t\le r,\ 1\le i\le d.
\EQNY
Next, for cases $a_2)$ and $a_3)$
\BQN\label{Eq: a23r}
\nonumber \lefteqn{\int_{\tilde A_i\cap\{z_{ij}>\lambda_i, z_{ik}\le\lambda_i\}}
\frac{\partial^2 g_h(\vk z)}{\partial z_{ij}\partial z_{ik}}\, d\vk z_i
=\int_{\tilde A_i\cap\{z_{ij}\le\lambda_i, z_{ik}>\lambda_i\}}
\frac{\partial^2 g_h(\vk z)}{\partial z_{ij}\partial z_{ik}}\, d\vk z_i } \\
&&=-
\int_{\{w_{i,r-1}''\le n\}}g_h(z_{ij}=\lambda_i, z_{ik}=\lambda_i)\, \frac{d\vk z_i}{dz_{ij} dz_{ik}}.
\EQN
Finally, for case $a_4)$
\BQNY
\int_{\tilde A_i\cap\{z_{ij}\le\lambda_i, z_{ik}\le\lambda_i\}}
\frac{\partial^2 g_h(\vk z)}{\partial z_{ij}\partial z_{ik}}\, d\vk z_i
=\int_{\{w_{ir}''\le n\}}g_h(z_{ij}=\lambda_i, z_{ik}=\lambda_i)\, \frac{d\vk z_i}{dz_{ij} dz_{ik}}.
\EQNY
This together with  \eqref{Eq: reexpp}--\eqref{Eq: a23r}  yields that
\BQN\label{Eq: Ere}
\tilde E_{ii}(j,k) &=&\int_{\tilde A_i\rq{}}\int_{\{w_{i,r-2}''\le n, w_{i,r-1}''=\IF\}} g_h(z_{ij}=\lambda_i, z_{ik}=\lambda_i)\, \frac{d\vk z}{dz_{ij} dz_{ik}}
\notag\\
&\quad&- \int_{\tilde A_i'}\int_{\{w_{i,r-1}''\le n, w_{ir}''=\IF\}}g_h(z_{ij}=\lambda_i, z_{ik}=\lambda_i)\, \frac{d\vk z}{dz_{ij} dz_{ik}} \notag\\
&=&\varphi(\lambda_i,\lambda_i;\delta_{ij,ik}^h) \notag\\
&& \times\left(\EE{\mathbb P\Big \{}\cap_{s=1, s\neq i}^d\{Z_{s(n-r+1)} >\lambda_s\}\cap(\vk Z_i''\in\{w_{i,r-2}''\le n, w_{i,r-1}''=\IF\})\Big\lvert \{Z_{ij}=\lambda_i, Z_{ik}=\lambda_i\}\Big\}\right. \notag
 \\&& -\left.  \EE{\mathbb P\Big \{}\cap_{s=1, s\neq i}^d\{Z_{s(n-r+1)} >\lambda_s\}\cap(\vk Z_i''\in\{w_{i,r-1}''\le n, w_{i,r}''=\IF\})\Big\lvert \{Z_{ij}=\lambda_i, Z_{ik}=\lambda_i\}\Big\}\right)
\EQN
establishing the {validity of  \eqref{Decomposition IIIET} for $i=l$}. 

\underline{b) Proof of \eqref{Decomposition IIIET} for $i<l$}. By $\tilde A_{il}''=\cap_{s=1,s\neq i,l}^d\cup_{t_1, \ldots, t_r=1}^n\{z_{st_1} >\lambda_s, \ldots, z_{st_r}>\lambda_s\}$ and $\tilde A_i$ in \eqref{Eq: reexpp}
\BQN\label{decomp. 1Dr}
\tilde E_{il}(j,k)=\int_{\tilde A_{il}''}\int_{\tilde A_i}\int_{\tilde A_l}
\frac{\partial^2 g_h(\vk z)}{\partial z_{ij}\partial z_{lk}}\, d\vk z.
\EQN
By decomposing the integral regions $\tilde A_i$ and $\tilde A_l$ according to $z_{ij}>,\le\lambda_i$ and $z_{lk} >,\le\lambda_l$, respectively, we obtain by similar arguments as for $E_{il}(j,k)$ that
\BQN\label{Eq: Ern}
\lefteqn{\tilde E_{il}(j,k) = \varphi(\lambda_i,\lambda_l;\delta_{ij,lk}^h)\EE{\mathbb P\Big \{}\cap_{s=1, s\neq i,l}^d\{Z_{s(n-r+1)} >\lambda_s\}\cap(\vk Z_i'\in\{w_{i,r-1}'\le n, w_{ir}'=\IF\})\notag
} \\&&\qquad\qquad \qquad\quad\  \cap(\vk Z_l'\in\{w_{l,r-1}'\le n, w_{lr}'=\IF\})\Big\lvert \{Z_{ij}=\lambda_i, Z_{lk}=\lambda_l\}\Big\},
\EQN
where  $w_{i1}' $ is introduced in \eqref{w'} and (similar notation for $w_{lt}'$ with respect to  $k$)
\BQNY
w_{it}'=\inf\{t_0\le n: z_{it_0}>\lambda_i, t_0\neq j, t_0>w_{i,t-1}'\},\quad 2\le t\le r,\ 1\le i\le d.
\EQNY
It follows then from \eqref{Eq: Ern}   that \eqref{Decomposition IIIET} holds.
 Consequently, the claim  in \eqref{ineq: abs} for $2<r\le n$ follows. 

Finally, in view of \eqref{cond: s-ind} we see that the \tm{indices} over the sum in \eqref{Decomposition II} and \eqref{Decomposition I(II)r}
are simplified to $1\le i< l\le d, 1\le j,k\le n$. Then the claim in \eqref{ineq: without abs} follows immediately from \eqref{Eq: En}, \eqref{Ineq: bivar norm} and \eqref{Eq: Ern}. This completes the proof.
\QED

   \prooftheo{T3} It is sufficient to present the proof of \eqref{Ineq: sharpest}.
In view of Lemma 4.2 in \cite{DebickiHJL14}, the claim in \eqref{Ineq: sharpest} for $r=1$ follows from condition \eqref{Assump.ind}. We shall
present next the proofs for a) $r=2$ and b) $2<r\le n$.  \\
\underline{a) Proof of \eqref{Ineq: sharpest} for $r=2$.}
It follows from  \eqref{Assump.ind}, \eqref{Decomposition I} and  \eqref{Decomposition II} that
\BQN\label{Eq: Decomposition}
\Delta_{(2)}(\vk u) =n \sum_{1\le i<l\le d}(\sigma_{il}^{(1)} -\sigma_{il}^{(0)})\int_0^1E_{il}\, dh,
\EQN
where $E_{il}:=E_{il}(1,1)$. Further, by \eqref{Assump.ind} and  \eqref{Eq: En} we have, with $\delta_{il}^h:=\delta_{i1,l1}^h$  (recall $\lambda_i:=-u_i, 1\le i\le d$)
\BQN\label{Ineq: En}
0\le \frac{E_{il}}{\varphi(-u_i, -u_l; \delta_{il}^h)}\le \cxL{\mathbb P\Big \{}\vk Z_i'\in\{w_{i1}'\le n, w_{i2}'=\IF\},\vk Z_l'\in\{w_{l1}'\le n, w_{l2}'=\IF\}\Big\}.
\EQN
Note that hereafter $w_{i1}', w_{i2}'$ and $w_{l1}', w_{l2}'$ are defined as in \eqref{w'} with respect to $j=k=1.$

Next, let $(\tilde Z_i, \tilde Z_l)$ be a bivariate standard normal random vector with correlation $\abs{\delta_{il}^h}$ and $u=\min_{1\le i\le {d}} u_i >0$.  It follows by Slepian's inequality in \cite{Slepian62}
and Lemma 2.3 in \cite{PicandsA} that
\BQNY
\cxL{\mathbb P\Big \{}Z_{ij}<-u_i, Z_{lk}<-u_l\Big\} &\le& \pk{\tilde Z_i<-u_i, \tilde Z_l<-u_l}\\
&\tm{\cxL{\le}}& \pk{-\tilde Z_i >u, -\tilde Z_l >u}
\le
\frac{(1+\abs{\delta_{il}^h})^2}{u^2}\varphi(u, u; \abs{\delta_{il}^h}), \quad j, k\le n,
\EQNY
implying thus
\BQNY
\lefteqn{\cxL{\mathbb P\Big \{}\vk Z_i'\in\{(w_{i1}',w_{i2}')=(2,\IF)\},\vk Z_l'\in\{(w_{l1}',w_{l2}')=(2,\IF)\}\Big\}}\\
&&=\cxL{\mathbb P\Big \{}Z_{i2}>-u_i, Z_{l2}>-u_i\Big\}\prod_{j=3}^n\cxL{\mathbb P\Big \{} Z_{ij}\le -u_i, Z_{lj}\le -u_l\Big\}  \le \Bigl(\frac{(1+\lvert \delta_{il}^h\lvert)^2}{u^2}\varphi(u,u;\lvert\delta_{il}^h\lvert)\Bigr)^{n-2}
\EQNY
and
\BQNY
\lefteqn{\cxL{\mathbb P\Big \{}\vk Z_i'\in\{(w_{i1}',w_{i2}')=(3,\IF)\},\vk Z_l'\in\{(w_{l1}',w_{l2}')=(2,\IF)\}\Big\} }\\
&&=\cxL{\mathbb P\Big \{}Z_{i2}<-u_i, Z_{l2}>-u_l, Z_{i3}>-u_i, Z_{l3}<-u_l\Big\}\prod_{j=4}^n \cxL{\mathbb P\Big \{}Z_{ij}<-u_i, Z_{lj}<-u_l\Big\}\\
&& \le \Bigl(\frac{(1+\lvert\delta_{il}^h\lvert)^2}{u^2}\varphi(u,u;\lvert\delta_{il}^h\lvert)\Bigl)^{n-2}.
\EQNY
Similarly, we may consider all $(n-1)^2$ cases  in \eqref{Ineq: En} for $w_{i1}'=w_{l1}'$  and $w_{i1}'\neq w_{l1}'$. Therefore, using further (4.6) in \cite{LiShao02} we have
\begin{align*}
E_{il}&\le (n-1)^2 \Bigl(\frac{(1+\lvert\delta_{il}^h\lvert)^2}{u^2}\varphi(u,u;\lvert\delta_{il}^h\lvert)\Bigl)^{n-2}\varphi(-u_i, -u_l;\delta_{il}^h)\\
&\le  \frac{(n-1)^2}{(2\pi)^{n-1}} u^{-2(n-2)} \frac{(1+\lvert\delta_{il}^h\lvert)^{2(n-2)}}{(1-\lvert\delta_{il}^h\lvert^2)^{(n-1)/2}}\expon{-\frac{(n-1)u^2}{1+\lvert\delta_{il}^h\lvert}}.
\end{align*}
Consequently, 
by \eqref{Eq: Decomposition} we have
\BQNY
\Delta_{(2)}(\vk u) &\le & n\sum_{1\le i<l\le d} (\sigma_{il}^{(1)} -\sigma_{il}^{(0)})_+\int_0^1E_{il}\, dh\notag \\
 & \le &\frac{n(n-1)^2}{(2\pi)^{n-1}}u^{-2(n-2)}\sum_{1\le i<l\le d} (\sigma^{(1)}_{il} -\sigma^{(0)}_{il})_+ \expon{-\frac{(n-1)u^2}{1+\rho_{il}}}\int_0^1 \frac{(1+\lvert\delta_{il}^h\lvert)^{2(n-2)} }{(1-\lvert\delta_{il}^h\lvert^2)^{(n-1)/2}}\, dh
 \notag\\ &=& \frac{n(n-1)^2}{(2\pi)^{n-1}} u^{-2(n-2)}\sum_{1\le i<l\le d} (A_{il}^{(2)})_+\expon{-\frac{(n-1)u^2}{1+\rho_{il}}}.
\EQNY
The last step follows since for $\delta_{il}^h=h(\sigma_{il}^{(1)} - \sigma_{il}^{(0)})+ \sigma_{il}^{(0)}$ we have
$\rho_{il}=\max(|\sigma_{il}^{(\cxL{0})}|,  |\sigma_{il}^{(1)}|)\ge\delta_{il}^h$ and
\BQN\label{Eq: integral(h)}
\int_0^1  \frac{(1+\lvert\delta_{il}^h\lvert)^{2(n-2)} }{(1-\lvert\delta_{il}^h\lvert)^2)^{(n-1)/2}}\, dh = \frac1{\sigma_{il}^{(1)} - \sigma_{il}^{(0)}}\int_{\sigma_{il}^{(0)}}^{\sigma_{il}^{(1)}}\frac{(1+\abs h)^{2(n-2)}}{(1-h^2)^{(n-1)/2}}\, dh.
 \EQN
\underline{b) Proof of \eqref{Ineq: sharpest} for $2<r\le n$.} By \eqref{Assump.ind} and \eqref{Decomposition I(II)r}
 \BQNY
\Delta_{(r)}(\u) = n\sum_{1\le i<l\le d}(\sigma_{il}^{(1)} -\sigma_{il}^{(0)})\int_0^1\tilde E_{il}\, dh,\EQNY
where $\tilde E_{il}:=\tilde E_{il}(1,1)$. Clearly, from \eqref{Eq: Ern} we have $\tilde E_{il}\ge0$. Further, similar arguments as for $E_{il}$ (consider the number of $w_{it}'=w_{ls}', s, t<r$) yield that
\BQNY
\frac{\tilde E_{il}}{\varphi(-u_i,-u_l; \delta_{il}^h)}&\le& \cxL{\mathbb P\Big\{}\vk Z_i'\in\{w_{i,r-1}'\le n, w_{ir}'=\IF\}, \vk Z_l'\in\{w_{l,r-1}'\le n, w_{lr}'=\IF\}\Big\}\\
&\le & (c_{n-1, r-1})^2\Bigl(\frac{(1+\lvert\delta_{il}^h\lvert)^2}{u^2}\varphi(u,u;\lvert\delta_{il}^h\lvert)\Bigl)^{n-r}.
\EQNY
Consequently, the claim in \eqref{Ineq: sharpest} for  $2<r\le n$ follows. 
{We complete the proof}.
 \QED

 \proofprop{T2} The lower bound follows directly from \netheo{T1}. Next we focus on the upper bound.
 We shall present  \cxL{below}  the proof for $r=2$. 
Hereafter, we adopt the same notation as in the proof of \netheo{T1}.
 %
Further, define 
\BQNY
  f(h)=  
  \exp\LT(\mathop{\sum_{1\le i<l\le d}}_{1\le j,k\le n} \frac1{H\cxL{\big(}(u_i+u_l)/2\big)} \mathcal C^h_{ij,lk}\RT),\ \ h\in[0, 1],
\EQNY
where $\mathcal C^h_{ij,lk}=\ln\fracl{\pi-2\arcsin(\sigma_{ij,lk}^{(0)})} {\pi-2\arcsin({\delta_{ij,lk}^h})}$ and 
$H(x)= \sqrt{2 \pi} e^{x^2/2}\E{(\mathcal{N}+x)_+}$ where $\mathcal{N}$ is a $N(0,1)$ random variable. 
It suffices to show that $ Q(\mathcal Z; \Gamma^h)  /f(h)$ is non-increasing in $h$, i.e., 
\BQN\label{Ineq: Qg}
 \frac{\partial Q(\mathcal Z; \Gamma^h) /\partial h}{ Q(\mathcal Z; \Gamma^h)} \le \frac{\partial f(h)/\partial h}{f(h)}, \quad h\in[0,1].
\EQN
Moreover, since
\BQN\label{Diff: g}
\frac{\partial f(h)/\partial h}{f(h)} 
= \mathop{\sum_{1\le i<l\le d}}_{1\le j,k\le n} \frac{2(\sigma_{ij,lk}^{(1)} - \sigma_{ij,lk}^{(0)})}{\cxL{\big(}\pi -2 \arcsin(\delta_{ij,lk}^h)\cxL{\big)\sqrt{1-(\delta_{ij,lk}^h)^2}}}\frac1{H\big((u_i+u_l)/2\big)}
\EQN
and, by \eqref{Decomposition II}
\BQN\label{Diff: -Q}
 \frac{\partial Q(\mathcal Z; \Gamma^h)}{\partial h} = 
 \mathop{\sum_{1\le i<l\le d}}_{1\le j,k\le n} (\sigma_{ij,lk}^{(1)} - \sigma_{ij,lk}^{(0)})E_{il}(j,k)\cxL.
\EQN
\cxL{Therefore, by \eqref{cond: s-ind},}  it is sufficient to show that
\BQN\label{Ineq: il}
E_{il}(j,k)\le \frac{2 Q(\mathcal Z; \Gamma^h)}{\big(\pi -2 \arcsin(\delta_{ij,lk}^h)\big)\sqrt{1-(\delta_{ij,lk}^h)^2}}\frac1{H\big((u_i+u_l)/2\big)},\quad  1\le i< l\le d, 1\le j,k\le n.
\EQN
From \eqref{Eq: En} we have (recall $\vk u=-\vk\lambda$)
\BQN\label{Ineq: 1}
 \frac{E_{il}(j,k)}{\varphi(u_i, u_l;  \delta_{ij,lk}^h)} &\le&  \mathbb P\left\{ \cap_{s=1, s\neq i,l}^d\{Z_{s(n-1)} >\lambda_s\}\cap (\vk Z_i'\in\{w_{i1}'\le n \}) \cap(\vk Z_l'\in\{w_{l1}'\le n \})\Big\lvert \{Z_{ij}=\lambda_i, Z_{lk} =\lambda_l\}\right\}\nonumber\\
&=&  \mathbb P\left\{\cap_{s=1, s\neq i,l}^d\{Z_{s(2)} {<} u_s\}\cap (\vk Z_i'\in\{v_{i1}'\le n \})  \cap(\vk Z_l'\in\{v_{l1}'\le n \})\Big\lvert \{Z_{ij}=u_i, Z_{lk} =u_l\}\right\},
\EQN
where $v_{i1}', v_{l1}'$ are defined by
\BQNY
v_{i1}'=\inf\{t: z_{it}<u_i, t\neq j\}, \quad  v_{l1}'=\inf\{t: z_{lt}<u_l, t\neq k\}.
\EQNY

Define next 
\BQNY
T_{ij}= \frac{(Z_{ij}- u_i) - \delta_{ij,lk}^h(Z_{lk} - u_l)}{1-(\delta_{ij,lk}^h)^2}, \quad T_{lk}= \frac{(Z_{lk}- u_l) - \delta_{ij,lk}^h(Z_{ij} - u_i)}{1-(\delta_{ij,lk}^h)^2}.
\EQNY
It follows  that the random vectors
$\vk Z_v^*=( Z_{vw}- \delta_{vw, ij}^h T_{ij} -\delta_{vw,lk}^h T_{lk}, 1\le w\le n), 1\le v  (\neq i,l)\le d,$  $\vk Z_i'^*=( Z_{it}- \delta_{it, ij}^h T_{ij} -\delta_{it,lk}^h T_{lk}, 1\le t\neq j\le n)$ and $\vk Z_l'^* =(Z_{lt}- \delta_{lt, ij}^h T_{ij} -\delta_{lt,lk}^h T_{lk}, 1\le t\neq k\le n)$
are independent of $(Z_{ij}, Z_{lk})$ and thus are independent of $(T_{ij}, T_{lk})$. Thus, 
by \eqref{Ineq: 1} and the fact that  $0\le \delta_{ij,lk}^h{<}1, (i-1)n+j <(l-1)n+k, \cxL{h\in[0,1]}$, \cxL{we have}
\BQN\label{Ineq: 2}
\lefteqn{ E_{il}(j,k)\frac{\pk{T_{ij}<0, T_{lk}<0}}{\varphi(u_i, u_l;  \delta_{ij,lk}^h)}}
\notag \\
& \le &  \pk{\cap_{s=1, s\neq i,l}^d\{Z_{s(2)}^* < u_s\}\cap (\vk Z_i'^*\in\{v_{i1}'\le n \}) \cap(\vk Z_l'^*\in\{v_{l1}'\le n \})\cap( T_{ij}<0)\cap( T_{lk}<0)}
\notag \\
& \le & \cxL{\mathbb P\Big\{}\cap_{s=1, s\neq i,l}^d\{Z_{s(2)}  < u_s\}\cap (\vk Z_i' \in\{v_{i1}'\le n \}) \cap(\vk Z_l' \in\{v_{l1}'\le n \})\cap(  Z_{ij}< u_i)\cap( Z_{lk}< u_l)\Big\} \notag \\
&{=}&
Q(\mathcal Z; \Gamma^h).
\EQN
Moreover, by Lemma 2.2 in \cite{Yan2009}
\BQNY
\frac{\pk{T_{ij}<0, T_{lk}<0}}{\varphi(u_i, u_l;  \delta_{ij,lk}^h)} \ge \frac{\pi -2\arcsin(\delta_{ij,lk}^h)}2 \sqrt{1-(\delta_{ij,lk}^h)^2} H\fracl{u_i+u_l}{2},
\EQNY
which together with \eqref{Ineq: 2} implies \eqref{Ineq: il}, hence the proof for $r=2$ is complete.\\
\cxL{For $2<r\le n$, we need to show that \eqref{Ineq: il} holds for $\tilde E_{il}(j,k)$.  This follows by similar arguments as for $r=2$, using the inequality \eqref{Eq: Ern} instead of \eqref{Eq: En}.}
\QED

\quad

\COM{
Since we can choose $T_d$ to bThe proof Next, we consider the closed interval $[0, T]$. For each $d\inn$, let $\mathcal T_d$ be a finite subset of $[0, T]$ such that $\mathcal T_d\subset \mathcal T_{d+1}$
and $\mathcal T_d \uparrow [0, T]$. By separability, we have $\cup_{d\ge 1}\mathcal T_d= [0, T]$ and thus ( $\toas$ denotes the convergence  almost surely)
\BQNY
\sup_{t\in \mathcal T_d} X_{r:n}(t)\toas \sup_{t\in [0, T]} X_{r:n}(t)
\EQNY
and
\BQNY
\sup_{t\in \mathcal T_d} Y_{r:n}(t) \toas \sup_{t\in [0, T]} Y_{r:n}(t)
\EQNY
hold as $d\to\IF$. Moreover,  the claim follows by \eqref{Ineq: discrete-cont} since the last convergence is monotone.}

\COM{\subsubsection{Proof of \nekorr{prop1}}
We  show that $p(x, T) = \ln \pk{\sup_{0\le t\le T} X_{r:n}(t) \le x}$ is sub-additive with respect to $T\ge0$.
\\
Given any $T, S\ge0$ let $\{\tilde X(t), t\in[0, T]\}\equivdis \{ X(t), t\in[0, T]\}$, independent of  $\{\tilde X(t+T), t\in[0, S]\}\equivdis \{ X(t), t\in[0, S]\}$. Note that $\E{X(0)X(t)} \ge0, t\ge0$ and $X$ is stationary. It follows from \neprop{Prop1} that
\BQNY
\lefteqn{\pk{\sup_{0\le t\le T+S} X_{r:n}(t) \le x} = \pk{\sup_{0\le t\le T} X_{r:n}(t) \le x, \sup_{T\le t\le T+S} X_{r:n}(t) \le x} }\\
&&\ge \pk{\sup_{0\le t\le T} \tilde X_{r:n}(t) \le x, \sup_{T\le t\le T+S} \tilde X_{r:n}(t) \le x} \\
&& = \pk{\sup_{0\le t\le T} \tilde X_{r:n}(t) \le x} \pk{\sup_{T\le t\le T+S} \tilde X_{r:n}(t) \le x} \\
&& = \pk{\sup_{0\le t\le T}  X_{r:n}(t) \le x} \pk{\sup_{0\le t\le S}  X_{r:n}(t) \le x}.
\EQNY
Thus, it follows that $p(x, T+S) \ge p(x, T)+p(x, S), T, S\ge 0$ which further implies that the limit $p(x)$ exists for all $x\inr$. The \wE{rest of the proof is the same} as for Proposition 3.1 in \cite{LiShao2004}.
\QED
}

 \prooftheo{ThmC}  
\EHH{First note that by \eqref{pr} 
 we have that} $\mathrm c_{r:n,\alpha}(x), x\inr$ defined by (with $Y_{r:n}(t):= X^*_{r:n}(t)+cZ^*(t)$)
$$
\mathrm c_{r:n,\alpha}(x)= -\lim_{T\to\IF}\frac1T\pk{\sup_{0\le t\le T} Y_{r:n}(t) \le x}, \quad x\in\R
$$
exists and \EHH{is} left-continuous.
Next, we show that $\mathrm c_{r:n,\alpha}(x)$ is right-continuous, which will be crucial for our proof. 
As in Theorem 3.1 (ii) in \cite{LiShao2005}, it suffices to  show that, for all $x\in\R, y>0, m\ge1, \theta\in(0,1)$ and $h\in(0, \IF)$
\BQN\label{r.cont.}
\pk{\sup_{0\le t\le m h} Y_{r:n}(t) \le x+y}\le \Phi^{{-m}}\fracl{-y+x\left(\sqrt{1+a_{h,\theta}^2}-1\right)}{a_{h,\theta}}\pk{\sup_{0\le t\le \theta mh} Y_{r:n}(t) \le x}.
\EQN
Let therefore $W_{k}, 1\le k\le m$ be independent $N(0,1)$ \cxL{random variables}  
which are further  independent of
the dual processes $X^*_i, Z^*,$ $1\le  i\le n$, and write, for simplicity, $a=a_{h,\theta}$. We have
\BQNY
{p_{h,\theta}(x, Y)}&:=&\pk{\max_{1\le k\le m}\sup_{(k-1)h\le t\le kh} \frac{Y_{r:n}(t,a W_k)}{\sqrt{1+a^2}} \le x} \\
& \ge &
\pk{\max_{1\le k\le m}\sup_{(k-1)h\le t\le kh} Y_{r:n}(t) \le x+y} \pk{\max_{1\le k\le m}aW_{k} \le -y+x\big(\sqrt{1+a^2}-1\big)} \\
&=&\pk{\max_{1\le k\le m}\sup_{(k-1)h\le t\le kh}Y_{r:n}(t) \le x+y} \Phi^{m}\fracl{-y+x\big(\sqrt{1+a^2}-1\big)}{a},
\EQNY
where {$\{Y_{r:n}(t,a  W_k), t\in [(k-1)h, kh]\}$ is the $r$th order statistics process generated by $\{Y_i(t)+a W_{k}, t\in [(k-1)h, kh)\}, {1\le i\le n}$.} Furthermore, 
it follows \cxL{by \eqref{continuity} and the monotonicity of $\rho(\cdot), \widetilde \rho(\cdot)$ that (set $I_k=[(k-1)h,kh)$)
\BQNY
&&\frac{\E{(Y(t)+aW_{ [t/h]+1})(Y(s)+aW_{[s/h]+1})}}{1+a^2} - \E{Y(\theta t)Y(\theta s)}\\
&&=\left\{\begin{array}{ll}
\frac{1-\rho(\abs{t-s})}{1+a^2}\left(a^2-\frac{\rho(\theta\abs{t-s}) -\rho(\abs{t-s})}{1-\rho(\abs{t-s})}\right) + c^2\left(\frac{\widetilde \rho(\abs{t-s})}{1+a^2} -\widetilde\rho(\theta\abs{t-s})\right),&t,s\in I_k;\\
 \frac{ \rho(\abs{t-s})}{1+a^2} - \rho(\theta\abs{t-s}) + c^2\left(\frac{\widetilde \rho(\abs{t-s})}{1+a^2} -\widetilde\rho(\theta\abs{t-s})\right), &  t\in I_k, s\in I_l, k\neq l\\
 \end{array}
 \right.\\
 &&\le 0, 
\EQNY
which implies by 
\neprop{Prop1} 
that}
 \BQNY
p_{h,\theta}(x,Y)
\le \pk{\max_{1\le k\le m}\sup_{(k-1)h\le t\le kh}Y_{r:n}(\theta t){\le x}}
\EQNY
establishing \eqref{r.cont.} and thus the continuity of $\mathrm c_{r:n,\alpha}(x)$ \cxL{follows}.
\EHH{In order} to complete the proof\cxL, it is suffices  to show that (set below $\Yr:=X_{r:n}(t)+cZ(t) $)
\BQN 
-\frac2{\alpha}\mathrm c_{r:n,\alpha}&\le& \liminf_{x\downarrow0}\frac{\ln \pk{\sup_{0\le t\le 1} \Yr \le x}}{\ln(1/x)}\notag \\
&\le& \label{sup.bound}
\limsup_{x\downarrow0}\frac{\ln \pk{\sup_{0\le t\le 1}  \Yr \le x}}{\ln(1/x)} \le -\frac2{\alpha}\mathrm c_{r:n,\alpha}.
\EQN
By the self-similarity \cxL{of the process $\widetilde Y$},  for any $x\in(0,1)$ we have 
\BQNY
\pk{\sup_{0\le t\le 2/\alpha \ln(1/x)} Y_{r:n}(t) \le 0}
&=& \pk{\sup_{x^{2/\alpha}\le t\le 1} \Yr \le 0}\\
& \le &  \pk{\sup_{x^{2/\alpha}\le t\le 1} \Yr \le x} \\
&\le& \frac{ \pk{\sup_{0< t\le 1} \Yr \le x} }{ \pk{\sup_{0 <t\le x^{2/\alpha}} \Yr \le x} }\\
&{=}& \frac{ \pk{\sup_{0< t\le 1} \Yr  \le x} }{ \pk{\sup_{0 <t\le1}\Yr \le 1} },
\EQNY
where the second inequality follows from   {\neprop{Prop1}} and the fact that {$\sigma_X(s,t)\ge0$}  {and $\sigma_Z(s,t)\ge0$ }
for all $s,t\ge0$.
Consequently, the lower bound in \eqref{sup.bound} follows since $\mathrm c_{r:n,\alpha}=\mathrm c_{r:n,\alpha}(0)$.
\cxL{Next,} we establish the upper bound in \eqref{sup.bound}. It follows that,  \cL{for $y>0$} \EHH{sufficiently small}
\BQNY
\lefteqn{ \frac1{(\alpha/2)h}\ln\pk{\sup_{0\le t\le h} Y_{r:n}(t) \le y} }\\
&=& \frac1{\alpha h/2}\ln\pk{\sup_{e^{-h}\le t\le 1} ( t^{-\alpha/2} \Yr ) \le y}
 \\
& \ge&  \frac1{\alpha h/2}\ln\pk{\sup_{e^{-h}\le t\le 1} \Yr  \le ye^{-\alpha h/2}}\\
& \ge& \frac1{\alpha h/2}\ln\pk{\sup_{0\le t\le 1} \Yr \le ye^{-\alpha h/2}} \\
&  = &{\frac{\alpha h/2 -\ln y}{\alpha h/2}}\frac1{\ln(1/(ye^{-\alpha h/2}))}\ln\pk{\sup_{0\le t\le 1} \Yr  \le ye^{-\alpha h/2}}.
\EQNY
Letting $h\to\IF$ in the above we obtain that
\BQNY
\limsup_{x\downarrow0}\frac{\ln \pk{\sup_{0\le t\le 1} \Yr \le x}}{\ln(1/x)} \le -\frac2{\alpha}\mathrm c_{r:n,\alpha}(y) \to -\frac2{\alpha}\mathrm c_{r:n,\alpha}, \quad y\downarrow0,
\EQNY
where the last step follows by the right-continuity of $\mathrm c_{r:n,\alpha}(x) $ at $0.$
Consequently,    \eqref{sup.bound} holds and thus the proof is complete.
\QED

\quad


 \section{Appendix}
\EHH{We present next the \cL{proof of} \eqref{Rem:1} and then present two lemmas which are used for the proof of \netheo{ThmB}.}
We conclude this section with the proof of \neprop{Prop1} and \netheo{ThmB}.  

 {\bf Proof of \eqref{Rem:1}.}
 The claim for $r=1$ follows from Theorem 2.1 in \cite{LuW2014}. For $2\le r\le n$, we see from the proof of \netheo{T1} that{,} it suffices to prove that $E_{ii}(j,k)\le0$ and $\tilde E_{ii}(j,k)\le0$ hold for all $1\le i\le d,1\le j<k\le n$. \\
{From Remark 2.5(3) in \cite{Lihaijun2009}, we see that all orthant tail dependence parameters of multivariate normal distributions are zero.} Therefore we have for instance for $j \neq1$ and $1\le i\le d$
\BQNY
\lefteqn{
\pk{\vk Z_{i}''\in\{w_{i1}''=\IF\}} - \pk{\vk Z_{i}''\in\{w_{i1}''=1, w_{i2}''=\IF\}}
}\\&&=
\bigl(1-2\pk{Z_{i1}>\lambda_i\lvert Z_{it}\le \lambda_i, t\neq1, j,k }\bigr)\pk{Z_{it}\le \lambda_i, t\neq1, j,k} \le0, \quad \lambda_i\to-\IF.
\EQNY
It follows then by \eqref{Eq: Ee}  that $E_{ii}(j,k)\le0$ for sufficiently large $u_i$ (equals $-\lambda_i)$.
Thus, we complete the proof for $r=2$.
Similar arguments show that $\tilde E_{ii}(j,k)\le0$  for  sufficiently large $u_i$   (recall \eqref{Eq: Ere}).
Consequently, the claim for $2<r\le n$ follows.
\cL{\QED}

For notational simplicity, we set $q=q(u)=u^{-2/\alpha}, u>0$ and write $[x]$ for the integer part of $x$.

\BL \label{LTD}
\EHH{Under the assumptions of \netheo{ThmB} \cxL{with $\gamma=0,$ }} 
then for any $a, T>0$
\BQN \label{cond: D}
 \limsup_{u\to\IF} \sum_{j =
[{T}/(aq)]}^{[\ve / \pk{X_{r:n}(0) >u}]} \pk{X_{r:n}(aqj)
>u \Big \lvert X_{r:n}(0) >u } \to 0, \quad \ve\downarrow 0.
\EQN
\EL
\prooflem{LTD}
By Lemma {2} in \cite{DebickiHJL14} (see the proof of {(3.20)} therein), for sufficiently large $u$
\BQNY
 {p_u(t)}:=\pk{X_{r:n}(t)>u \Big \lvert X_{r:n}(0)>u} \le  2\pk{X_{r:r}(t)>u, X_{r:r}(0)>u \Big \lvert X_{r:r}(0)>u}.
 \EQNY
Since further $X(t) - \rho(t)X(0)$ is independent of $X(0)$,  we have for some constant $K>0$ {(the value of $K$ might change below from line to line)}
 \BQN
 p_u(t)&\cxL{\le}& 2^{r+1}\left(  \pk{X(t)>X(0)>u \Big \lvert X(0)>u} \right)^r\notag\\
 &\le&  2^{r+1}\left(  \pk{X(t) - \rho(t)X(0) >u (1-\rho(t)), X(0) >u \Big \lvert X(0)>u} \right)^r \notag\\
 &=& 2^{r+1} \left( 1-\Phi\left(u\sqrt{\frac{1-\rho(t)}{1+\rho(t)}}\right)\right)^r
 \notag \\
&\le& K u^{-r} \fracl{1-|\rho(t)|}{1+|\rho(t)|}^{-r/2} \expon{-\frac{ru^2}{2}\frac{1-|\rho(t)|}{1+|\rho(t)|}},
\label{eqn: D}
 \EQN
\cxL{the} last inequality follows by the Mill\rq{}s ratio inequality $1-\Phi(x)\le 1/(\sqrt{2\pi}x)\expon{-x^2/2}, x>0$.

Now we choose a function  $ g=g(u)$ such that $\lim_{u\to\IF}g(u)=\IF,  {|\rho(g(u))|} = u^{-2}$. Further it follows from  $u^{-2}\ln g(u) = o(1)$ that $g(u) \le \exp(\epsilon' u^2)$ for some $0<\epsilon'<r/2(1-|\rho(T)|)/(1+|\rho(T)|)$
(recall that $|\rho(T)|<1$; see \cxL{\cite{leadbetter1983extremes}, p.\,86}) 
 and sufficiently large $u$.
Next, we split the sum in \eqref{cond: D} at $aqj = g(u)$. The first term {is}
\BQNY
\lefteqn{
\sum_{j =[T/(aq)]}^{[g(u)/ (aq)]} \pk{X_{r:n}(aqj) >u \Big \lvert X_{r:n}(0) >u }}\\
&& \le  K\frac{g(u)}{aq} u^{-r} \fracl{1-|\rho(T)|}{1+|\rho(T)|}^{-r/2} \expon{-\frac{ru^2}{2}\frac{1-|\rho(T)|}{1+|\rho(T)|}}
 \\
& &\le {K  u^{2/\alpha -r}\expon{\epsilon' u^2-\frac{ru^2}{2}\frac{1-|\rho(T)|}{1+|\rho(T)|}}}\to 0, \quad u\to\IF.
\EQNY
For the remaining term, by Lemma \cL1 in \cite{DebickiHJL14}
\BQNY
\lefteqn{
\sum_{j = [g(u)/ (aq)]}^{[\ve /\pk{X_{r:n}(0) >u}]} \pk{X_{r:n}(aqj) >u \Big \lvert X_{r:n}(0) >u }}\\
&&\le  K \frac{\ve}{\pk{X_{r:n}(0) >u}} u^{-r} \fracl{1-u^{-2}}{1+u^{-2}}^{-r/2} \expon{-\frac{ru^2}{2}\frac{1-u^{-2}}{1+u^{-2}}} \\
&& \le  K \ve \expon{-\frac{ru^2}{2}\left(\frac{1-u^{-2}}{1+u^{-2}}-1\right)}\\
&& \le K \ve, \quad u\to\IF.
\EQNY
Therefore, the claim follows by taking $\ve\downarrow 0$.
\QED

\EHH{Next,} with the notation as in \eqref{Asym: finite-time} we set
\BQN
\label{def: T}
T= T(u) = \frac {1}{ c_{n,r}\mathcal A_{r,\alpha}} (2 \pi)^{\frac r2}u^{r-\frac 2\alpha} \expon{ \frac{r u^2}2},\ \ u>0.
\EQN
\BL \label{LTD2}
Let $T=T(u)$ be defined as in \eqref{def: T} and $a>0, 0 < \lambda < 1$ be any given constants.
\EHH{Under the assumptions of \nelem{LTD}}
for any $0\le s_1 < \cdots < s_p < t_1 < \cdots < t_{p'}$ in
 $ \{ aqj: j\in \EE{\mathbb{Z}}, 0\le aqj \le T\} $ with $t_1 - s_p \ge \lambda T$
 \BQN\label{Asymp.ind}
&&\Big\lvert
 \pk{ \cap_{ i=1}^p\{X_{r:n}(s_i) \le u\}, \cap_{ j=1}^{p'}\{X_{r:n}(t_j) \le u\} }\notag \\
&&\quad  -
 \cxL{\mathbb P\Big\{}\cap_{ i=1}^p\{X_{r:n}(s_i) \le u\}\Big\} \pk{\cap_{ j=1}^{p'}\{X_{r:n}(t_j) \le u\}}
\Big\lvert    \to0, \quad u\to\IF.
 \EQN
\EL

\prooflem{LTD2}
Denote
\BQNY
X_{ij}= X_j(s_i)\mathbb I\{i\le p\} +  X_j(t_{i-p})\mathbb I\{p<i\le p+p'\}, \quad 1\le i\le p+p',\ 1\le j\le n,
\EQNY
and $\{Y_{ij}, 1\le i\le p, 1\le j\le n\}\equivdis \{X_{ij}, 1\le i\le p, 1\le j\le n\}$, independent of
$\{Y_{ij}, p+1\le i\le p+p', 1\le j\le n\}\equivdis \{X_{ij}, p+1\le i\le p+p', 1\le j\le n\}$.
Applying \netheo{T3} with $X_{i(n-r+1)}= X_{r:n}(s_i)\mathbb I\{i\le p\} + X_{r:n}(t_{i-p})\mathbb I\{p<i\le p+p'\}$ {and  $Y_{i(n-r+1)}= Y_{r:n}(s_i)\mathbb I\{i\le p\} + Y_{r:n}(t_{i-p})\mathbb I\{p<i\le p+p'\}$},  it follows by similar arguments {as for} Lemma 8.2.4 in \cite{leadbetter1983extremes} that, the left-hand side of \eqref{Asymp.ind} is bounded from above by
\BQNY
\lefteqn{
K {u^{-2(r-1)}}\fracl Tq\sum_{\lambda T\le t_j-s_i \le T} \expon{-\frac{{r u^2}}{1+\abs{\rho(t_j-s_i)}}}\int_0^{\abs{\rho(t_j-s_i)}}\frac{(1+\abs h)^{2(r-1)}}{(1-h^2)^{r/2}}\, dh
} \\
&& \le Ku^{-2(r-1)}\fracl Tq \sum_{\lambda T \le aqj \le T} \abs{\rho(aqj)} \expon{-\frac{ru^2}{1+\abs{\rho(aqj)} }}\quad {\rm for\ } u\ {\rm large},
\EQNY
where $K$ is some constant. The rest of the \aE{proof} consists of the same arguments as that of
Lemma 12.3.1 in \cite{leadbetter1983extremes} using further the following asymptotic relation (recall \eqref{def: T})
 \BQNY  \label{Asym: T}
 u^2  \cL=
 \frac2r\ln T + \left(\frac 2{r\alpha}-1\right)\ln \ln T +
 \ln\left(\fracl r2^{1-2/(r\alpha)} \frac {(c_{n,r}\mathcal A_{r,\alpha})^{2/r}}{2\pi}\right) {(1+o(1))}, \quad \EHH{u\to \IF},
\EQNY
hence the  proof is complete. \QED
\COM{ (the constant $K$ below may be different from line to line). \\
Next, letting $\gamma(t) = \sup\{\abs{\rho(s)}\ln s: s\ge t\}, t\ge 1$,  we have that $\abs{\rho(t)} \le \gamma(t) /\ln t$ and $\gamma(t) \le M$ for some positive constant $M$ and all sufficiently large $t$.
Therefore, by \eqref{Asym: T}
\BQNY
\expon{{-\frac{ru^2}{1+\abs{\rho(aqj)}}}}
& \le&
\expon{-ru^2\left(1 - \frac{\gamma(\lambda T)}{\ln(\lambda T)}\right)}
\\
&\le&  { K \expon{-ru^2}\le K T^{-2}(\ln T)^{r-2/\alpha}}
\EQNY
holds for $T$ large.
Consequently,
\BQNY
\lefteqn{
 K  u^{-2(r-1)}\frac Tq\sum_{\lambda T \le aqj \le T} \abs{\rho(aqj)} \expon{{-\frac{ru^2}{1+\abs{\rho(aqj)}}}}
 } \\
 && \le  K u^{-2(r-1)}\fracl Tq^2\frac1{T/q}\sum_{\lambda T \le aqj \le T} \abs{\rho(aqj)} {\ln(aqj) }\frac1{\ln(\lambda T)}T^{-2}(\ln T)^{r-2/\alpha} \\
&& \le K \frac1{T/q}\sum_{\lambda T \le aqj \le T} \abs{\rho(aqj)} \ln(aqj),
\EQNY
which tends to 0 as $T\to\IF$ since $\rho(t)\ln t = o(1)$. }

\EHH{Below $W$ is \cxL{an} $N(0,1)$ random variable independent of any other random element \cxL{involved}. }

\proofprop{Prop1}
{We shall \cxL{first} present the proof of \eqref{Slep:Order1} 
for any finite set $\mathcal T_d$ containing $d$ elements such that $\mathcal T_d\subset [0, T]$. We write $\mathcal{T}_d =\{t_1,\ldots,t_d\}\subset[0, T]$.
Further we define  $f(t_i)=\sqrt{\sigma_X(t_i, t_i)+c^2\sigma_Z(t_i, t_i)}=\sqrt{\sigma_Y(t_i, t_i)+c^2\sigma_Z(t_i, t_i)}$ and
\BQNY
X^*_{ij} := \frac{X_j(t_i)+cZ(t_i)}{f(t_i)} , \quad Y_{ij}^* := \frac{Y_j(t_i)+cZ(t_i)}{f(t_i)}, \quad 1 \leq i \leq d, \ 1 \leq j \leq n.
\EQNY
 Then $X^*_{ij}$ and $Y^*_{ij}$ are $N(0,1)$ distributed, and 
\BQNY
\pk{\sup_{t\in\mathcal T_d}\big( X_{n-r+1:n}(t)+cZ(t)\big) >u} = \pk{\cup_{i=1}^d\{X^*_{i (r)}>u_i\}},\quad u_i:= \frac u{ f(t_i)}, \ 1\le i\le d.
\EQNY
Noting \EE{that} $\{Z(t),t\ge0\}$ is  independent of  $\{X(t),t\ge0\}$ and  $\{Y(t),t\ge0\}$ we have
\BQNY
\E{X_{ij}^*X^*_{ik}} &=& \frac{\E{X_j(t_i)X_k(t_i) + c^2Z^2(t_i)}}{ (f(t_i))^2}\\
&=& \frac{{\sigma_X(t_i, t_i)}\mathbb I\{j=k\}+c^2 \sigma_Z(t_i, t_i)}{(f(t_i))^2}= \E{Y^*_{ij}Y^*_{ik}},\ 1\le i\le d, \ 1\le j,  k\le n
\EQNY
and
\BQNY
\E{X_{ij}^*X_{lk}^*}&=&
 \frac{\sigma_X(t_i, t_l) \mathbb I\{j=k\} +c^2 \sigma_Z(t_i, t_l)}{ f(t_i) f(t_l)}\\
 &\le&  \frac{\sigma_Y(t_i, t_l)  \mathbb I\{j=k\} +c^2 \sigma_Z(t_i, t_l)}{ f(t_i) f(t_l)}=\E{Y_{ij}^*Y_{lk}^*}, \quad 1\le i<l\le d,\ 1\le  j,  k\le n.
\EQNY
Therefore, by   \eqref{ineq: discrete}
\BQNY\label{Ineq: discrete-cont}
\pk{\sup_{t\in\mathcal T_d} \cxL{\big(}Y_{r:n}(t) +cZ(t)\big)>u} \le \pk{\sup_{t\in\mathcal T_d} \big( X_{r:n}(t) +cZ(t)\big)>u}.
\EQNY
The passage from $\mathcal T_d$ to $[0,T]$ is standard and therefore we omit the details.  \cxL{We thus complete the proof of \eqref{Slep:Order1}. \\
Next, for \eqref{Slep:Order2}, we denote instead 
 $f(t_i)=\sqrt{\sigma_Z(t_i, t_i)+c^2\sigma_X(t_i, t_i)}=\sqrt{\sigma_Y(t_i, t_i)+c^2\sigma_Y(t_i, t_i)}$ and
\BQNY
X^*_{ij} := \frac{Z_j(t_i)+cX(t_i)}{f(t_i)} , \quad Y_{ij}^* := \frac{Z_j(t_i)+cX(t_i)}{f(t_i)}, \quad 1 \leq i \leq d, \ 1 \leq j \leq n.
\EQNY
Then the \eW{rest of the proof} is the same as that for \eqref{Slep:Order1}.}  \QED

\prooftheo{ThmB} a)
In view of {Theorem 10 in \cite{Albin1990}, since \eqref{Asym: finite-time} and Lemmas \ref{LTD} and \ref{LTD2} hold for the $r$th order statistics process
$\{X_{r:n}(t),t\ge0\}$, we have  for $T=T(u)$ defined as in \eqref{def: T}
\BQNY
 \lim_{u\to\IF}\pk{\sup_{t\in[0, T(u)]}X_{r:n}(t) \le u + \frac x{ru} } = \expon{-e^{-x}}, \quad x\in\R.
\EQNY
 Expressing $u$ in term{s} of $T$ using \eqref{def: T}  we obtain
the required claim  for any $x\in\R$, with $a_{r, T}, b_{r, T}$  given as in
\eqref{Normalized constant}; the uniform convergence in $x$ follows since all functions (with respect to $x$) are
continuous, bounded and increasing.

\def\cTT{\EHH{c_T}}
b) The proof follows from the main arguments of Theorem 3.1 in \cite{MittalY1975} by showing {that}, for any $\ve>0$ and $x\in\R$
\BQN\label{Ineq:IF}
\Phi(x-\ve)&\le& \liminf_{T\to\IF}  \pk{\EHH{M_X(T)}\le  \cTT b_{r, T} +\sqrt{\rho(T)} x} \notag\\
&\le& \limsup_{T\to\IF} \pk{M_X(T) \le\cTT b_{r, T} +\sqrt{\rho(T)} x} \le \Phi(x+\ve),
\EQN
\EHH{where $M_X(t):=\sup_{ t\in [0, T]} X_{r:n}(t)$ and $\cTT:=\sqrt{1-\rho(T)}$.}
We start with the proof of the first inequality. Let {$\rho^*(t), t\ge0$ be a correlation function of a stationary Gaussian process such that } $\rho^*(t)= 1-2\abs{t}^\alpha + o(\abs{t}^\alpha)$ as $t\to0$. Then there exists some $t_0>0$ such that for $T$ large
\BQN\label{Ineq:cons1}
\rho^*(t)\EE{c_T^2} +\rho(T) \le \rho(t),\quad 0\le t\le t_0.
\EQN
Denote by $\{Y_k(t), t\ge0\}, k\in  \cxL{\mathbb{N}}$ 
independent centered stationary Gaussian process\cxL{es}  with a.s. continuous sample paths and common covariance function $\rho^*(\cdot)$, and define $\{Y(t), t\ge0\}$ by
 \BQN\label{Constr.GP1}
 Y(t)=\sum_{k=\cxL1}^\IF Y_k(t)\mathbb{I}\{t\in [\cxL{(k-1)}t_0,  \cxL{k}t_0)\},\ \ t\ge0.
\EQN
 It follows from \eqref{Ineq:cons1} that for $T$ \EHH{sufficiently} large
\BQNY
\E{X(s)X(t)} \ge \E{\big(\EE{c_T}Y(s)+\sqrt{\rho(T)}W\big)\big( \EE{c_T} Y(t)+\sqrt{\rho(T)}W\big)}, \quad s, t\ge0.
\EQNY
Therefore, by \neprop{Prop1}
\BQNY
\pk{\EHH{M_X(T)}\le\cTT b_{r, T} +\sqrt{\rho(T)} x} &\ge & \pk{\cTT \EHH{M_Y(T)} +\sqrt{\rho(T)}W\le\cTT b_{r, T}+\sqrt{\rho(T)} x} \\
 &\ge& \Phi(x-\ve) \left( \pk{\sup_{ t\in [0, t_0]} Y_{r:n}(t) \le b_{r, T}+\ve\sqrt{\rho(T)} }\right)^{[T/t_0]+1}.
\EQNY
Noting that $a=\inf_{0<t\le t_0} (1-\rho^*(t))\abs t^\alpha >0$, we have by  Theorem 1.1 in \cite{DebickiHJL14} (see \cxL{also} \eqref{Asym: finite-time})
\BQNY
\lim_{T\to\IF}\frac{\pk{\sup_{ t\in [0, t_0]} Y_{r:n}(t) >b_{r, T} +\ve\sqrt{\rho(T)}}}
{t_0c_{n,r} b^{2/\alpha}_{r, T}\big(1-\Phi(b_{r, T} +\ve\sqrt{\rho(T)})\big)^r }
=   2^{1/\alpha}\mathcal A_{r,\alpha}. 
\EQNY
Consequently, since $\gamma=\IF$ we have
\BQNY
\lefteqn{\lim_{T\to\IF}([T/t_0]+1) \ln \pk{\sup_{ t\in [0, t_0]} Y_{r:n}(t) \le b_{r, T}+\ve\sqrt{\rho(T)}}
}\\
&&=-\lim_{T\to\IF}\frac T{t_0} \pk{\sup_{ t\in [0, t_0]} Y_{r:n}(t) \tm{>} b_{r, T}+\ve\sqrt{\rho(T)}}\\
&&= -\lim_{T\to\IF} T c_{n,r} 2^{1/\alpha}\mathcal A_{r,\alpha}b^{2/\alpha}_{r, T}\big(1-\Phi(b_{r, T} +\ve\sqrt{\rho(T)})\big)^r \\
&&= 0 
\EQNY
establishing the first inequality in \eqref{Ineq:IF}.\\
Next, we consider the last inequality in \eqref{Ineq:IF}. Note that, \cxL{by} the convexity of $\rho(\cdot)$, there is a separable stationary Gaussian process $\{Y(t), t\in[0, T]\}$ with correlation function given by
\BQN\label{def:rho\rq{}}
\tilde \rho(t)= \frac{\rho(t)-\rho(T)}{1-\rho(T)},\quad t\in[0, T].
\EQN
We have
\BQNY
\EHH{M_X(T)} {=}\cTT \EHH{M_Y(T)} + \sqrt {\rho(T)}W.
\EQNY
Therefore
\BQN\label{FactorN}
\pk{\EHH{M_X(T)}\le\cTT  b_{r, T} +\sqrt{\rho(T)} x}
& =&\int_{-\IF}^\IF\pk{\EHH{M_Y(T)}\le b_{r, T} +\frac{\sqrt{\rho(T)}}{\cTT } (x-u)}\varphi(u)\,du \notag \\
&\le& \Phi(x+\ve) + \pk{\EHH{M_Y(T)}\le b_{r, T} -\ve\frac{\sqrt{\rho(T)}}{\cTT }},
\EQN
which means that we only need to prove
\BQNY
\lim_{T\to\IF} \pk{\EHH{M_Y(T)}\le b_{r, T} -\ve\sqrt{\rho(T)}} =0. 
\EQNY
To this end, using again the convexity of $\tilde\rho\cxL{(\cdot)}$, we construct a separable stationary Gaussian process $\{Z(t), t\in[0,T]\}$ with the  correlation function \cxL{(recall $\tilde\rho(\cdot)$ in \eqref{def:rho\rq{}})}
\BQN\label{def:sigma}
\sigma(t)= \max\Big(\tilde \rho(t), \tilde \rho\big(T\exp\big(-\sqrt{\ln T}\big)\big) \Big),\quad t\in[0, T].
\EQN
Again by \neprop{Prop1}, we have
\BQN\label{Ineq:63}
\pk{\EHH{M_Y(T)}\le b_{r, T} -\ve\sqrt{\rho(T)}}\le \pk{\EHH{M_Z(T)}\le b_{r, T} -\ve\sqrt{\rho(T)}}.
\EQN
Now we make a grid as follows. {Let}  $I_1, \ldots, I_{[T]}$  be $[T]$ consecutive unit intervals with an interval of length $\delta$ removed from the right-hand side of each one with $\delta\in(0,1)$ given, \cxL{and}
$$\mathcal G_T=\big\{k(2\ln T)^{-3/\alpha},\ k\cxL{\inn}\big\}\cap \big(\cup_{i=1}^{[T]}  I_i\big).$$ %
 {It follows} from  Theorem 10 in \cite{Albin1990} and Theorem 1.1 in \cite{DebickiHJL14} that, $\sup_{ t\in [0, T]} Z_{r:n}(t)$ and $\sup_{ t\in \mathcal G_T} Z_{r:n}(t)$ have the same asymptotic distribution and thus we only need to show that
\BQNY
\lim_{T\to\IF}\pk{\sup_{ t\in  \mathcal G_T} Z_{r:n}(t)\le b_{r, T} -\ve\sqrt{\rho(T)}}=0. 
\EQNY
Let {$\{Z'_{r:n}(t),t\ge0\}$} be generated by $\{Z'(t), t\in[0,T]\}$ which is again a separable stationary process with the correlation function \cxL{(recall $\sigma(\cdot)$ in \eqref{def:sigma})}
\BQNY
\sigma'(t)=\frac{\sigma(t)-\sigma(T)}{1-\sigma(T)}, \quad t\in[0, T].
\EQNY
\EHH{Analogous to the derivation of \eqref{FactorN} we obtain }
\BQNY
\lefteqn{\pk{\sup_{ t\in \mathcal G_T} Z_{r:n}(t)\le b_{r, T} -\ve\sqrt{\rho(T)}}}\\
&&=\pk{ \sqrt{1-\sigma(T)} \max_{t\in \mathcal G_T}Z_{r:n}'(t) + \sqrt{\sigma(T)}\cL{W} \le b_{r, T} -\ve\sqrt{\rho(T)} } \\
&&\le \Phi\left(-\frac12\ve \fracl{\rho(T)}{\sigma(T)}^{1/2}\right) +\pk{\max_{t\in \mathcal G_T}Z'_{r:n}(t) \le b_{r,T}+\frac{b_{r,T}\sigma(T)}{\sqrt{1-\sigma(T)}(1+\sqrt{1-\sigma(T)})}-\frac{\ve\sqrt{\rho(T)}}{2\sqrt{1-\sigma(T)}}},
\EQNY
which tends to 0 as $T\to\IF$. The proof of it is the same as that of Theorem 3.1 in \cite{MittalY1975}, by using instead Theorem 1.1 in \cite{DebickiHJL14} and our \netheo{T3}. {Consequently, the last inequality in \eqref{Ineq:IF} follows by \eqref{FactorN} and  \eqref{Ineq:63}.} We complete the proof for $\gamma=\IF$.

 c)  Given $\delta\in(0,1)$, take $I_1, \ldots, I_{[T]}$  as in b). For $\{Y_k(t), t\ge0\}, k\inn$ independent copies of $X$ define 
 \BQNY
 Y(t):= \sum_{k=1}^\IF Y_k(t)\mathbb I\{t\in[k-1,k)\}, \quad t\ge0
 \EQNY
and 
\BQNY
X_*(t):=\sqrt{1-\rho_*(T)} Y(t) +\sqrt{\rho_*(T)}W,\quad t\in \cup_{k=\cxL{1}}^{[T]}{I_k},
\EQNY
where $\rho_*(T)= \gamma/\ln T$. 
The rest of the proof is similar to that as for Theorem 2.1 in \cite{ZTEH12} \cxL{by using our \netheo{T3} instead of Berman\rq{}s inequality. We omit the details.} \\
Combining all the arguments for the three cases above, we complete the proof of \netheo{ThmB}. \QED


\section*{Acknowledgments} 
Supported in part from SNSF 200021-140633/1 and the project RARE -318984 (an FP7 Marie Curie IRSES Fellowship).
The first author also acknowledges  partial support by NCN Grant No 2013/09/B/ST1/01778 (2014-2016).

\bibliographystyle{plain}

 \bibliography{BermanOrderST}
\end{document}